\documentclass[12pt]{amsart}
\usepackage{amsmath}
\usepackage{amsxtra}
\usepackage{amscd}
\usepackage{amsthm}
\usepackage{amsfonts}
\usepackage{amssymb}
\usepackage{eucal}

\newtheorem{thm}{Theorem}[section]
\newtheorem{lem}{Lemma}[section]
\newtheorem{cor}{Corollary}[section]
\newtheorem{prop}{Proposition}[section]

\theoremstyle{definition}

\theoremstyle{remark}
\newtheorem{rem}{Remark}[section]

\numberwithin{equation}{section}

\begin{document}

\newcommand{\thmref}[1]{Theorem~\ref{#1}}
\newcommand{\secref}[1]{Section~\ref{#1}}
\newcommand{\lemref}[1]{Lemma~\ref{#1}}
\newcommand{\propref}[1]{Proposition~\ref{#1}}
\newcommand{\corref}[1]{Corollary~\ref{#1}}
\newcommand{\remref}[1]{Remark~\ref{#1}}
\newcommand{\eqnref}[1]{(\ref{#1})}
\newcommand{\exref}[1]{Example~\ref{#1}}

\newcommand{\nc}{\newcommand}
\nc{\on}{\operatorname} \nc{\Z}{{\mathbb Z}} \nc{\C}{{\mathbb C}}
\nc{\R}{{\mathbb R}} \nc{\boldD}{{\mathbb D}} \nc{\oo}{{\mf O}}
\nc{\N}{{\mathbb N}} \nc{\bib}{\bibitem} \nc{\pa}{\partial}
\nc{\F}{{\mf F}} \nc{\rarr}{\rightarrow}
\nc{\larr}{\longrightarrow} \nc{\al}{\alpha} \nc{\ri}{\rangle}
\nc{\lef}{\langle} \nc{\W}{{\mc W}} \nc{\gam}{\ol{\gamma}}
\nc{\Q}{\ol{Q}} \nc{\q}{\widetilde{Q}} \nc{\la}{\lambda}
\nc{\ep}{\epsilon} \nc{\g}{\mf g} \nc{\h}{\mf h} \nc{\n}{\mf n}
\nc{\bb}{\mf b} \nc{\A}{{\mf a}} \nc{\G}{{\mf g}} \nc{\D}{\mc D}
\nc{\Li}{{\mc L}} \nc{\La}{\Lambda} \nc{\is}{{\mathbf i}}
\nc{\V}{\mf V} \nc{\bi}{\bibitem} \nc{\NS}{\mf N}
\nc{\dt}{\mathord{\hbox{${\frac{d}{d t}}$}}} \nc{\E}{\mc E}
\nc{\ba}{\tilde{\pa}} \nc{\half}{\frac{1}{2}}
\def\smapdown#1{\big\downarrow\rlap{$\vcenter{\hbox{$\scriptstyle#1$}}$}}
\nc{\mc}{\mathcal} \nc{\mf}{\mathfrak} \nc{\ol}{\fracline}
\nc{\el}{\ell} \nc{\etabf}{{\bf \eta}} \nc{\zetabf}{{\bf
\zeta}}\nc{\x}{{\bf x}} \nc{\xibf}{{\bf \xi}} \nc{\y}{{\bf y}}
\nc{\WW}{\mc W} \nc{\SW}{\mc S \mc W} \nc{\sd}{\mc S \mc D}
\nc{\hsd}{\widehat{\mc S\mc D}} \nc{\parth}{\partial_{\theta}}
\nc{\cwo}{\C[w]^{(1)}} \nc{\cwe}{\C[w]^{(0)}}
\nc{\hf}{\frac{1}{2}} \nc{\hsdzero}{{}^0\widehat{\sd}}
\nc{\hsdpp}{{}^{++}\widehat{\sd}}
\nc{\hsdpm}{{}^{+-}\widehat{\sd}}
\nc{\hsdmp}{{}^{-+}\widehat{\sd}}
\nc{\hsdmm}{{}^{--}\widehat{\sd}} \nc{\gltwo}{{\rm
gl}_{\infty|\infty}} \nc{\btwo}{{B}_{\infty|\infty}}
\nc{\htwo}{{\h}_{\infty|\infty}} \nc{\hglone}{\widehat{\rm
gl}_{\infty}} \nc{\hgltwo}{\widehat{\rm gl}_{\infty|\infty}}
\nc{\hbtwo}{\hat{B}_{\infty|\infty}}
\nc{\hhtwo}{\hat{\h}_{\infty|\infty}} \nc{\glone}{{\rm gl}_\infty}
\nc{\gl}{{\rm gl}} \nc{\ospd}{\mc B} \nc{\hospd}{\widehat{\mc B}}
\nc{\pd}{\mc P} \nc{\hpd}{\widehat{\pd}} \nc{\co}{\mc O}
\nc{\Oe}{\co^{(0)}} \nc{\Oo}{\co^{(1)}} \nc{\sdzero}{{}^0{\sd}}
\nc{\hz}{\hf+\Z} \nc{\vac}{|0 \rangle} \nc{\K}{\mf k}
\nc{\bhf}{\bf\hf}

\advance\headheight by 2pt

\nc{\fb}{{\mathfrak b}} \nc{\fg}{{\mathfrak g}}
\nc{\fh}{{\mathfrak h}}  \nc{\fk}{{\mathfrak k}}
\nc{\fl}{{\mathfrak l}} \nc{\fn}{{\mathfrak n}}
\nc{\fp}{{\mathfrak p}}

\nc{\fu}{u} \nc{\cF}{{\mathcal F}}

\title[Character Formula for general linear superalgebra]
{Character Formula for \\ Infinite Dimensional Unitarizable
Modules\\ of the General Linear Superalgebra}

\author[Shun-Jen Cheng]{Shun-Jen Cheng$^1$}
\thanks{$^1$Partially supported by NSC-grant 91-2115-M-002-007 of the R.O.C.}
\address{Department of Mathematics, National Taiwan University, Taipei,
Taiwan 106} \email{chengsj@math.ntu.edu.tw}

\author[Ngau Lam]{Ngau Lam$^{2}$}
\thanks{$^{2}$Partially supported by NSC-grant 91-2115-M-006-013 of the
R.O.C.}
\address{Department of Mathematics, National Cheng Kung University, Tainan,
Taiwan 701} \email{nlam@mail.ncku.edu.tw}

\author[R.~B.~Zhang]{R.~B.~Zhang$^{3}$}
\thanks{$^{3}$Supported by the Australian Research Council.}
\address{School of Mathematics and Statistics, University of Sydney,
New South Wales 2006, Australia} \email{rzhang@maths.usyd.edu.au}

\begin{abstract} The Fock space of $m+p$ bosonic and $n+q$
fermionic quantum oscillators forms a unitarizable module of the
general linear superalgebra $gl_{m+p|n+q}$. Its tensor powers
decompose into direct sums of infinite dimensional irreducible
highest weight $gl_{m+p|n+q}$-modules.  We obtain an explicit
decomposition of any tensor power of this Fock space into
irreducibles, and develop a character formula for the irreducible
$gl_{m+p|n+q}$-modules arising in this way.

\vspace{.3cm}

\noindent{\bf Key words:} Lie superalgebra, unitarizable
representations, Howe duality, character formula.

\vspace{.3cm}

\noindent{\bf Mathematics Subject Classifications (1991)}: 17B67.

\end{abstract}
\maketitle

\section{Introduction}
The Fock space of $m+p$ bosonic and $n+q$ fermionic quantum
oscillators (see Subsection \ref{comments} for definition) with
the standard inner product furnishes a unitarizable complex
representation of the real form $\fu(m, p|n, q)$ of the general
linear superalgebra $gl_{m+p|n+q}$.  This representation
decomposes into a direct sum of infinite dimensional irreducible
representations which are of highest weight type with respect to
an appropriate choice of a Borel subalgebra. Because of the
unitarity, any tensor power of the representation is also
semi-simple with all irreducible sub-representations being
unitarizable highest weight representations. We shall characterize
the irreducible sub-representations and determine their structure.

In recent years there have been considerable activities (see,
e.g., \cite{D} for references) in the physics community to study
unitarizable highest weight representations of Lie superalgebras.
This is motivated by applications of such representations in
quantum field theory. For example, the symmetry algebra of the yet
largely conjectural $M$-theory is closely related to
$osp_{1|32}(\R)$ \cite{T}. An understanding of the unitarizable
highest weight representations of this Lie superalgebra will help
to solve mysteries of $M$-theory. It has also been recognized
\cite{DFLV} that some real forms of simple basic classical Lie
superalgebras provide the conformal superalgebras of higher
dimensional space-time manifolds with extended supersymmetries.
The unitarizable highest weight representations of these Lie
superalgebras thus describe the spectra of possible elementary
particles existing in such space-times.

The problem of determining the possible unitarizable irreducible
highest weight representations of real forms of simple Lie
superalgebras was investigated by a number of people (see \cite{D}
and references therein), with  the most systematical study given
in \cite{J}. However, a classification analogous to the
Enright-Howe-Wallach \cite{EHW} classification of unitarizable
positive energy irreducible representations for ordinary real
simple Lie algebras has yet to be achieved (see Subsection
\ref{comments}).

A demanding but physically more important problem is to understand
the structure of the unitarizable irreducible representations.
Recall that a character formula for the unitarizable irreducible
highest weight representations of real forms of simple Lie
algebras \cite{EHW} was given in \cite{E} some fifteen year ago.
In an earlier publication \cite{CZ}, two of the authors studied
the irreducible representations arising from the decomposition of
the tensor powers of the oscillator representations of the
orthosymplectic superalgebras. By using results of \cite{DES, E},
a character formula for these irreducible representations was
derived. In this paper we investigate the case of the general
linear superalgebra.

It is known from \cite{H1} that $\fu(d)$ and $\fu(m, p|n, q)$ form
a dual reductive pair on the $d$-th tensor power of the Fock space
of $m+p$ bosonic and $n+q$ fermionic quantum oscillators. We
explore the duality between the complexifications of these Lie
(super)algebras to obtain in \thmref{hw+*} an explicit
decomposition of the tensor power of the Fock space into
irreducible $gl_d\times gl_{m+p|n+q}$-modules. The Howe duality
again as in \cite{CL, CZ} forms the key ingredient and further
enables us to compute the characters for the irreducible
$gl_{m+p|n+q}$-representations of \thmref{hw+*}. This result is
presented in \thmref{chWLa}.  Another application of the Howe
duality is the computation of the tensor product decomposition of
any two such unitarizable modules, which is the contents of
\thmref{tensor1}.

Here is an outline of the paper.  In \secref{finite-duality} we
discuss the $(gl_d, gl_{m|n})$-dualities on
$S(\C^d\otimes\C^{m|n})$ and its graded dual space. The material
is largely known, but the highest weight vectors in the graded
dual of $S(\C^d\otimes\C^{m|n})$ given in \lemref{glpq-hwv} have
not been computed previously as far as we are aware of. In
\secref{infinite-duality} we study the
$gl_{m+p|n+q}$-representations furnished by tensor powers of the
Fock space of $m+p$ bosonic and $n+q$ fermionic quantum
oscillators. In Subsection \ref{unitarity} we show that such
representations are unitarizable and their irreducible
sub-representations are infinite dimensional highest weight
representations, and in Subsection \ref{glglduality**} we obtain
the explicit decomposition of the $d$-th tensor power of the Fock
space with respect to the semi-simple multiplicity free action of
$gl_d\times gl_{m+p|n+q}$. \secref{branching} gives the
$gl_{m+p|n+q}\rightarrow gl_{p|q}\times gl_{m|n}$ branching rule
for the infinite dimensional unitarizable irreducible
$gl_{m+p|n+q}$-representations arising from the decomposition of
tensor powers of the Fock space. In \secref{character-formula} we
develop a character formula for these infinite dimensional
irreducible $gl_{m+p|n+q}$-representations in terms of hook Schur
functions. Finally in \secref{tensorproduct} we calculate the
tensor product decomposition of two such irreducible
$gl_{m+p|n+q}$-modules that appear in our decompositions of tensor
powers.

\section{Tensorial representations of
general linear superalgebra}\label{finite-duality}
This section presents some results on the $(gl_d, gl_{m|n})$-dualities
on $S(\C^d\otimes\C^{m|n})$ and its graded dual vector space. The material
contained here will be important for the remainder of the paper.
\subsection{Preliminaries}\label{Preliminaries}
We work on the field $\C$ of complex numbers throughout the paper.
Let $gl_d$ denote the Lie algebra of all complex $d\times d$
matrices. Let $\{ e^1,\ldots,e^d\}$ be the standard basis for
$\C^d$. Denote by $e_{ij}$ the elementary matrix with $1$ in the
$i$-th row and $j$-th column and $0$ elsewhere. Then
$\fh_d=\sum_{i=1}^d\C e_{ii}$ is a Cartan subalgebra, while
$\fb_d=\sum_{1\le i\le j\le d}\C e_{ij}$ is the standard Borel
subalgebra containing $\h_d$. The weight of $e^i$ is denoted by
$\tilde{\epsilon}_i$ for $1\le i\le d$.

Let $\C^{m|n}=\C^{m|0}\oplus\C^{0|n}$ denote the $m|n$-dimensional
superspace. The superspace of complex linear transformations on
$\C^{m|n}$ has a natural structure of a Lie superalgebra \cite{K},
which we will denote by $gl_{m|n}$. Choose a basis $\{e_1,\cdots,e_m\}$
for the even subspace $\C^{m|0}$ and a basis $\{f_1,\cdots,f_n\}$
for the odd subspace $\C^{0|n}$, then $\{e_1,\cdots,e_m, f_1,\cdots,f_n\}$
is a homogeneous basis for $\C^{m|n}$. We may regard $gl_{m|n}$ as
consisting of $(m+n)\times(m+n)$ matrices relative to this basis.
Denote by $E_{ij}$ the
elementary matrix with $1$ in the $i$-th row and $j$-th column and
$0$ elsewhere. Then $\fh_{m|n}=\sum_{i=1}^{m+n}\C E_{ii}$ is a Cartan
subalgebra, while $\fb_{m|n}=\sum_{1\le i\le j\le m+n}\C E_{ij}$ is the
standard Borel subalgebra containing $\h_{m|n}$. We shall denote
the weights of $e_i$ and $f_j$ by $\epsilon_i$ and $\delta_j$
respectively for $i=1, \cdots, m$, and $j=1, \cdots, n$.

By a partition $\la$ of length $k$ we mean a non-increasing finite
sequence of non-negative integers $(\la_1,\cdots,\la_k )$. We will
let $\la'$ denote the transpose of the partition $\la$. For
example, if $\la=(4,3,1,0,0)$, then the length of $\la$ is $5$ and
$\la'=(3,2,2,1)$. By a generalized partition of length $k$, we
shall mean a non-increasing finite sequence of integers
$(\la_1,\cdots,\la_k )$.  In particular, every partition is a
generalized partition of non-negative integers.  Corresponding to
each generalized partition $\la=(\la_1,\dots,\la_d)$, we will
define $\la^*:=(-\la_d,\dots,-\la_1)$. Then $\la^*$ is also a
generalized partition.

We regard a finite sequence $\la=(\la_1,\cdots,\la_d )$ of complex
numbers as an element of the dual vector space $\fh^*_d$ of
$\fh_d$ defined by $\la(e_{ii})=\la_i$, for $i=1,\cdots,d$. Denote
by $V^\la_d$ the irreducible $gl_d$-module with highest weight
$\la$ relative to the standard Borel subalgebra $\fb_d$.
Similarly, we shall also regard a finite sequence of complex
numbers $\la=(\la_1,\cdots,\la_{m+n} )$ as an element of the dual
vector space $\fh^*_{m|n}$ of $\fh_{m|n}$ such that $\lambda(E_{j
j})=\la_j$, $1\le j\le m+n$.  We denote by $V^\la_{m|n}$ the
irreducible $gl_{m|n}$-module with highest weight $\la$ relative
to the standard Borel subalgebra $\fb_{m+n}$.

\subsection{The $(gl_d,  gl_{m|n})$-duality on $S(\C^d\otimes\C^{m|n})$}

Recall that the natural action of the Lie superalgebra $gl_d\times
gl_{m|n}$ on $\C^{d}\otimes\C^{m|n}$ induces an action on the
supersymmetric tensor algebra $S(\C^{d}\otimes\C^{m|n})$. This
action is completely reducible and multiplicity free \cite{H1, S,
CW1, CW2}. Indeed the pair $(gl_d,gl_{m|n})$ forms a dual
reductive pair on $S(\C^{d}\otimes\C^{m|n})$ in the sense of Howe
\cite{H1, H2}.
\begin{thm}\cite{CW1}\label{hw-glmn}
Under the $gl_d\times gl_{m|n}$-action, $S(\C^d\otimes\C^{m|n})$
decomposes into
\begin{equation}\label{glgl-duality}
S(\C^d\otimes\C^{m|n})\cong\sum_{\la} V^\la_d\otimes
V^{\widetilde\la}_{m|n},
\end{equation}
where the sum in \eqref{glgl-duality} is over all partitions $\la$
of length $d$ subject to the condition $\la_{m+1}\le n$, and
\begin{equation}\label{wtglmn}
\widetilde\la=(\la_1, \cdots, \la_m; <\la'_1-m>, \cdots,
<\la'_n-m>).
\end{equation}
Here $\la'$ is the transpose partition of $\la$, and $<r>$ stands
for $r$, if $r\in\N$, and $0$ otherwise.
\end{thm}
\begin{rem}
The condition $\la_{m+1}\le n$ is considered to be automatically
satisfied by every generalized partition $\la$ of length $d$ if
$m\ge d$.
\end{rem}

We shall need an explicit formula for the joint highest weight
vectors of the irreducible $gl_d\times gl_{m|n}$-module
$V^\la_d\otimes V^{\widetilde\la}_{m|n}$ inside
$S(\C\otimes\C^{m|n})$. (See also \cite{OP} and \cite{N} for
different descriptions of these vectors.) We set
\begin{equation}\label{generators}
x_l^i:=e^i\otimes e_l,\quad \eta_k^i:=e^i\otimes f_k,
\end{equation}
for $i=1, \cdots, d$, \ $l=1, \cdots, m$, and $k=1, \dots, n$.  We
will denote by $\C[\x,\etabf]$ the polynomial superalgebra
generated by \eqnref{generators}. By identifying
$S(\C^{d}\otimes\C^{m|n})$ with $\C[\x,\etabf]$ the commuting pair
$(gl_d, gl_{m|n})$ can be realized in terms of the following first
order differential operators ($1\le i,i'\le d$, $1\le s,s'\le m$
and $1\le k,k'\le n$):
\begin{align}
\phi(e_{i i'}):=&\sum_{j=1}^m x_{j}^{i}\frac{\partial}{\partial
x_j^{i'}}+\sum_{j=1}^n\eta_j^{i}\frac{\partial}{\partial\eta_j^{i'}},
\label{glpq} \allowdisplaybreaks\\
\phi(E_{s s'}):=&\sum_{j=1}^d x_{s}^{j}\frac{\partial}{\partial
x_{s'}^j},\quad\quad \phi(E_{m+k, m+k'}):=
\sum_{j=1}^d\eta_{k}^{j}\frac{\partial}{\partial\eta_{k'}^j},\label{glmn1}\\
\phi(E_{s, m+ k}):=& \sum_{j=1}^d
x_s^j\frac{\partial}{\partial\eta_k^j}, \quad\quad \phi(E_{m+k,
s}):= \sum_{j=1}^d\eta_k^j\frac{\partial}{\partial
x_s^j}.\nonumber
\end{align}
Straightforward calculations show that $\phi(e_{i j})$, $1\le i,
j\le d$, and $\phi(E_{a b})$, $1\le a, b \le m+n$, satisfy the
same commutation relations as the elementary matrices $e_{i j}$
and $E_{a b}$ respectively. Thus \eqnref{glpq} spans a copy of
$gl_d$, and \eqnref{glmn1} a copy of $gl_{m|n}$.

For $1\le r\le {\rm min}(d,m)$, we define
\begin{equation}\label{deltar}
\Delta_r:={\rm det}\begin{pmatrix}
x_{1}^{1}&x_{2}^{1}&\cdots&x_{r}^{1}\\
x_{1}^{2}&x_{2}^{2}&\cdots&x_{r}^{2}\\
\vdots&\vdots&\vdots&\vdots\\
x_{1}^{r}&x_{2}^{r}&\cdots&x_{r}^{r}\\
\end{pmatrix}.
\end{equation}
If $d>m$, we consider the following determinant of an
$r\times r$ matrix for every  $m< r\le d$:
\begin{equation} \label{eqdet}
\Delta_{k,r}:={\rm det}
\begin{pmatrix}
x_1^1&x_1^2&\cdots &x_1^r\\
x_2^1&x_2^2&\cdots &x_2^r\\
\vdots&\vdots&\cdots &\vdots\\
x_m^1&x_m^2&\cdots &x_m^r\\
\eta_k^1&\eta_k^2&\cdots &\eta_k^r\\
\eta_k^1&\eta_k^2&\cdots &\eta_k^r\\
\vdots&\vdots&\cdots &\vdots\\
\eta_k^1&\eta_k^2&\cdots &\eta_k^r\\
\end{pmatrix},\quad k=1,\ldots,n.
\end{equation}
That is, the first $m$ rows are filled by the vectors $(x_j^1,
\ldots, x_j^r)$, for $j=1, \ldots, m$, in increasing order and the
last $r-m$ rows are filled with the same vector
$(\eta_k^1,\ldots,\eta_k^r)$. Here the determinant of an $r\times
r$ matrix
$$
A :=\begin{pmatrix} a_1^1 & a_1^2 &\cdots&a_1^r\\
a_2^1&a_2^2&\cdots&a_2^r\\
\vdots&\vdots&\cdots &\vdots\\
a_r^1&a_r^2&\cdots&a_r^r\\
\end{pmatrix},
$$
with matrix entries possibly involving Grassmann variables, is by
definition the expression $\sum_{\sigma\in
S_r}(-1)^{l(\sigma)}a_1^{\sigma(1) }a_2^{\sigma (2) }\cdots
a_r^{\sigma (r) }$, where $l(\sigma)$ is the length of $\sigma$ in
the symmetric group $S_r$.

Observe that both $\Delta_r$ and $\Delta_{k,r}$ (if defined) are
weight vectors of $gl_d\times gl_{m|n}$. Their $gl_d$-weights are
respectively
\begin{align}\label{wtd}\begin{array}{l l}
wt_d(\Delta_r)=(1, \cdots, \underbrace{1}_r, 0, \cdots,
0),\\
wt_d(\Delta_{k,r})=(1, \cdots, \underbrace{1}_r, 0, \cdots, 0),
\end{array}\end{align}
while the $gl_{m|n}$-weights are respectively
\begin{align}\label{wtmn}\begin{array}{l l}
wt_{m|n}(\Delta_r)=(1, \cdots, \underbrace{1}_r, 0, \cdots,
0),\\
wt_{m|n}(\Delta_{k,r})=(1, \cdots, \underbrace{1}_m, 0, \cdots,0,
\underbrace{r-m}_{m+k}, 0,\cdots, 0),
\end{array}\end{align}
Corresponding to each partition $\la$ of length $d$ satisfying the
condition $\la_{m+1}\le n$, we define
\begin{eqnarray}\label{hwvform}
\Delta_\la
&:=&\left\{\begin{array}{l l}\Delta_{\la'_1}\Delta_{\la'_2}
\cdots\Delta_{\la'_{\la_1}},
&\quad \mbox{if}\ \la_1^\prime\le m,\\
\prod_{k=1}^{\la_{m+1}}\Delta_{k,\la'_k}\prod_{j=1+{\la_{m+1}}}^{\la_1}\Delta_{\la'_j},
&\quad  \mbox{if}\  \la_1^\prime> m.
\end{array}\right.
\end{eqnarray}
\begin{lem}\cite{CW1}\label{aux2}\label{glpmn-duality}
The space of $gl_d \times gl_{m|n}$ highest weight vectors in the
submodule $V^\la_d\otimes V^{\widetilde\la}_{m|n}$ of $ \C[\x,
\etabf]$ is $\C \Delta_\la$.
\end{lem}

\subsection{The $(gl_d, gl_{p|q})$-duality on
$S({\C^d}^*\otimes{\C^{p|q}}^*)$}\label{glglduality*} Let us
denote by ${\C^{p|q}}^*$ the dual of the natural $gl_{p|q}$-module
$\C^{p|q}$, and by ${\C^d}^*$ the dual of the natural
$gl_d$-module $\C^d$. Then the $gl_d\times gl_{p|q}$-action on
${\C^d}^*\otimes{\C^{p|q}}^*$ induces a $gl_d\times
gl_{p|q}$-action on $S({\C^d}^*\otimes{\C^{p|q}}^*)$

If $S^k(W)$ denotes the set of all homogeneous elements of degree
$k$ in the supersymmetric tensor algebra of the superspace $W$,
then $S^k({\C^d}^*\otimes{\C^{p|q}}^*) \cong
S^k({\C^d}\otimes\C^{p|q})^*$, and thus
$S({\C^d}^*\otimes{\C^{p|q}}^*) \cong\sum_k
S^k({\C^d}\otimes\C^{p|q})^*$. Therefore it follows from the
decomposition \eqref{glgl-duality} that the $gl_d\times
gl_{p|q}$-action on $S({\C^d}^*\otimes{\C^{p|q}}^*)$ is also
semi-simple and multiplicity free. Furthermore, we have the
following decomposition $
S({\C^d}^*\otimes{\C^{p|q}}^*)\cong\sum_{\la}(V^\la_d)^*\otimes
(V^{\widetilde\la}_{p|q})^*,$ where $\la$ is summed over all
partitions of length $d$ subject to the condition $\la_{p+1}\le
q$. Clearly, $(V^\la_d)^*\cong V_d^{\la^*}$. Also,
$(V^{\widetilde\la}_{p|q})^*\cong V^{\widehat{\la^*}}_{p|q}$,
where $\widehat{\la^*}$ is the negative of the lowest weight of
$V^{\widetilde\la}_{p|q}$.  We shall give an explicit formula for
$\widehat{\la^*}$ in equation \eqref{hatla}.

Since the supertrace ${\rm Str}$ is trivial on the derived algebra
of $gl_{p|q}$, one may twist any action of $gl_{p|q}$ by any
scalar multiple of the supertrace. This is to say that, if $X\in
gl_{p|q}$ acts on a space, then we may define a new action of $X$
on this space by $X+\gamma{\rm Str}(X)$ instead, where
$\gamma\in\C$. This in particular allows us to twist the standard
action of $gl_{p|q}$ on $S(\C^{d*}\otimes\C^{p|q*})$ by $-d{\rm
Str}$. Under this twisted action of $gl_{p|q}$ the space
$S({\C^d}^*\otimes{\C^{p|q}}^*)$ decomposes into
\begin{equation}\label{glgl*}
S({\C^d}^*\otimes{\C^{p|q}}^*)\cong\sum_{\la}V^\la_d\otimes
V^{-d{\bf 1}+\widehat\la}_{p|q},
\end{equation}
where ${\bf 1}:=(1, \cdots, \underbrace{1}_p, -1, \cdots,
-1)$. Here the summation in $\la$ is over all generalized partitions of
non-positive integers with length $d$ subject to $\la_{d-p}\ge
-q$. Observe that $\la^*$ is a partition of length $d$ satisfying $(\la^*)_{p+1}\le q$.
\begin{rem} \label{pged} For any generalized partition $\la$ of length $d$,
the condition $\la_{d-p}\ge -q$ is considered to be automatically satisfied if
$p\ge d$.
\end{rem}
\begin{rem}
Hereafter we shall always refer to this twisted action of
$gl_{p|q}$ when considering $S({\C^d}^*\otimes{\C^{p|q}}^*)$ as a
$gl_d\times gl_{p|q}$-module.
\end{rem}

Let ${e^1},\ldots,{e^d}$ be the standard basis of $\C^d$.  Let
${e^1}^*,\ldots,{e^d}^*$ be a basis for the contragredient
$gl_d$-module ${\C^d}^*$. We require the two bases to be dual in the sense that
${e^i}^*(e^j)=\delta_{i j}$ for all $i, j \in\{1, \cdots d\}$.
Similarly, we let $e_1,\ldots,e_p,$ $f_1,\ldots,f_q$ denote the
standard homogeneous basis for the natural $gl_{p|q}$-module
$\C^{p|q}$ and $e_1^*, \ldots, e_p^*,$ $f_1^*,\ldots,f_q^*$ denote
the dual basis for the contragredient $gl_{p|q}$-module
${\C^{p|q}}^*$. For $1\le l\le d$, \ $1\le i\le p$ and $1\le j\le
q$, we set
\begin{equation}\label{generators*}
y_i^l:={e^l}^*\otimes e_i^*,\quad \zeta_j^l:={e^l}^*\otimes f_j^*,
\end{equation}
which form a basis for ${\C^d}^*\otimes{\C^{p|q}}^*$.

We will denote by $\C[\y,\zetabf]$ the polynomial superalgebra
generated by \eqref{generators*}. Let $e_{ij}$, $1\le i,j\le d$
and $E_{a b}$, $1\le a,b\le p+q$ be the bases respectively for
$gl_d$ and $gl_{p|q}$ consisting of elementary matrices. Then the
action of the commuting pair $(gl_d,gl_{p|q})$ on $\C[\y,\zetabf]$
can be realized in terms of first order differential operators as
follows ($1\le i,j\le d$, $1\le r,r'\le p$ and $1\le s,s'\le q$):
\begin{align}
\bar\phi(e_{ij}):=&-\sum_{k=1}^p y_{k}^{j}\frac{\partial}{\partial
y_k^{i}}-\sum_{k=1}^q\zeta_k^{j}\frac{\partial}{\partial\zeta_k^{i}};
\label{barphid}\\
\bar\phi(E_{rr'}):=&-\sum_{l=1}^d \frac{\partial}{\partial
y_{r}^l}y_{r'}^{l},\quad\quad
\bar\phi(E_{s+p,s'+p}):=\sum_{l=1}^d\frac{\partial}{\partial\zeta_s^l}
\zeta_{s'}^{l}, \label{barphi} \\
\bar\phi(E_{s+p,r}):=&-\sum_{l=1}^d
\frac{\partial}{\partial\zeta_s^l} y_r^l,\quad\quad
\bar\phi(E_{r,s+p}):=\sum_{l=1}^d\frac{\partial}{\partial y_r^l}
\zeta_s^l.\nonumber
\end{align}
It is straightforward to show that the $\bar\phi(e_{i j})$ and
$\bar\phi(E_{a b})$ satisfy the same commutation relations as the
$e_{i j}$ and $E_{a b}$ respectively. Furthermore, the elements of
\eqref{barphid} commute with those of \eqref{barphi}.

For $1\le r\le {\rm min}(d,p)$, we define the following
determinant of an $r\times r$ matrix:
\begin{equation}\label{deltar*}
{\Delta}_r^*:={\rm det}\begin{pmatrix}
&y_{p}^{d}&y_{p-1}^{d}&\cdots&y_{p-r+1}^{d}\\
&y_{p}^{d-1}&y_{p-1}^{d-1}&\cdots&y_{p-r+1}^{d-1}\\
&\vdots&\vdots&\vdots&\vdots\\
&y_{p}^{d-r+1}&y_{p-1}^{d-r+1}&\cdots&y_{p-r+1}^{d-r+1}\\
\end{pmatrix}.
\end{equation}
For $1\le r\le d$, we define
\begin{equation} \label{eqdet*}
{\Delta}_{k,r}^*:=
\zeta_k^d\zeta_k^{d-1}\cdots \zeta_k^{d-r+1}\\
\quad k=1,\ldots,q.
\end{equation}

It is clear that the $\Delta^*_r$ and $\Delta^*_{k, r}$ are all
$gl_d$ highest weight vectors with respect to the standard Borel
subalgebra $\fb_d$. They are also weight vectors under the action
of $gl_d\times gl_{p|q}$ with the $gl_d$-weights respectively
given by
\begin{align}\label{gldwts}
\begin{array}{r c l}
wt_d(\Delta^*_r)&=&(0, \cdots, 0, \underbrace{-1}_{d+1-r}, \cdots, -1),\\
wt_d(\Delta_{k, r})&=&(0, \cdots, 0, \underbrace{-1}_{d+1-r},
\cdots, -1),
\end{array}\end{align}
and the $gl_{p|q}$-weights respectively given by
\begin{align}\label{glpqwts}
\begin{array}{r c l}
wt_{p|q}(\Delta^*_r)&=&-d {\bf 1} + (0, \cdots, 0, \underbrace{-1}_{p+1-r},
\cdots,  \underbrace{-1}_{p}, 0, \cdots, 0),\\
wt_{p|q}(\Delta_{k, r})&=&-d {\bf 1} + (0, \cdots, 0,
\underbrace{-r}_{p+k}, 0, \cdots, 0).
\end{array}\end{align}

Let $\la=(\la_1, \cdots, \la_d)$ be a generalized partition of
non-positive integers subject to the condition $\la_{d-p}\ge -q$.
Then $\mu:=\la^*$ is a partition satisfying the condition
$\mu_{p+1}\le q$. We let $\mu^\prime$ denote the transpose
partition of $\mu$.  Define
\begin{eqnarray}
\Delta^*_\la&:=& \left\{\begin{array}{l l}
\prod_{k=1}^{\mu_1}\Delta^*_{q+1-k, \mu^\prime_k}, & \mbox{if}\
\mu_1\le q, \\
\prod_{k=1}^{q}\Delta^*_{q+1-k, \mu^\prime_k}
\prod_{l=q+1}^{\mu_1}\Delta^*_{\mu^\prime_l}, & \mbox{if}\ \mu_1
>q.
\end{array}\right.\label{Deltala}
\end{eqnarray}
\begin{lem}\label{glpq-hwv}
Let $\la=(\la_1, \cdots, \la_d)$ be a generalized partition of
non-positive integers subject to the condition $\la_{d-p}\ge -q$.
Then $\Delta^*_\la$ is a non-zero highest weight vector in $\C[\y,
\zetabf]$ with respect to the joint actions of $gl_d$ and
$gl_{p|q}$. The $gl_d$-weight of  $\Delta^*_\la$ is $\la$, while
the $gl_{p|q}$-weight of $\Delta^*_\la$ is $-d {\bf 1}+\hat\la$ with
\begin{eqnarray} \label{hatla} \hat{\la}&:=&- \left(\langle
\mu_p-q\rangle, \langle \mu_{p-1}-q\rangle, \cdots,\langle
\mu_1-q\rangle,  \mu^\prime_q, \mu^\prime_{q-1}, \cdots,
\mu^\prime_1\right), \end{eqnarray} where $\mu=\la^*$.
\end{lem}
\begin{proof} Note that $\bar\phi(e_{i j})$, for all $1\le i, j\le d$,
act on $\C[\y, \zetabf]$ by derivations. Thus the product of any
subset of the $gl_d$ highest weight vectors $\Delta^*_r$,
$\Delta^*_{k, r}$, $r=1, 2, \cdots, min(d, p)$, \  $k=1, 2, \cdots,
q$,  is also a $gl_d$ highest weight vector. Hence so is
$\Delta^*_\la$.

Obviously $\Delta^*_\la$ is a highest weight vector with respect
to the action of the subalgebra $gl_p\times gl_q$  of $gl_{p|q}$.
Consider the action of $\bar\phi(E_{p, p+1})
=\sum_{l=1}^d\frac{\partial}{\partial y_p^l} \zeta_1^l$ on
$\Delta^*_\la$. When $-\la_{d}=\mu_1\le q$, $\Delta^*_\la$ does
not involve any of the variables $y^l_p$, thus  $\bar\phi(E_{p,
p+1})\Delta^*_\la=0$. By using the equation
\begin{equation*}
\Delta^*_{1, r} \sum_{l=1}^d \zeta_1^l \frac{\partial}{\partial
y_p^l}\Delta_{s}=0,  \quad r\ge s,
\end{equation*}
we also easily show  that $\Delta^*_\la$ is annihilated by
$\bar\phi(E_{p, p+1})$ when $\mu_1> q.$ This proves that
$\Delta^*_\la$ is indeed a $gl_d\times gl_{p|q}$ highest
weight vector. The rest of the lemma
easily follows from equations \eqref{gldwts} and \eqref{glpqwts}.
\end{proof}

To summarize this subsection, we combine \lemref{glpq-hwv} with
equation \eqref{glgl-duality} into the following theorem.

\begin{thm}\label{hw*} Under the $gl_d \times gl_{p|q}$-action,
$\C[\y, \zetabf]$ decompose into
\begin{equation*}
\C[\y, \zetabf]\cong\sum_{\la} V^\la_d\otimes V^{-d {\bf 1}+
\widehat\la}_{p|q},
\end{equation*}
where $\la$ is summed over all generalized partitions of
non-positive integers with length $d$ subject to $\la_{d-p}\ge
-q$. The space of highest weight vectors in $V^\la_d\otimes V^{-d
{\bf 1}+ \widehat\la}_{p|q}$ is given by $\C\Delta_{\la}^*$.
\end{thm}

\section{The $(gl_d,gl_{m+p|n+q})$-Duality on
$S({\C^d}\otimes\C^{m|n}\oplus{\C^d}^*\otimes{\C^{p|q}}^*)$}\label{infinite-duality}
\subsection{The $gl_d\times gl_{m+p|n+q}$-action on
$S({\C^d}\otimes\C^{m|n}\oplus{\C^d}^*\otimes{\C^{p|q}}^*)$} We
described the semi-simple multiplicity free actions of $gl_d\times
gl_{m|n}$  on $S({\C^d}\otimes\C^{m|n})$ and $gl_d\times gl_{p|q}$
on $S({\C^d}^*\otimes{\C^{p|q}}^*)$ in the last section. Through
the obvious isomorphism between $S({\C^d}\otimes\C^{m|n})
\otimes_{\C}S({\C^d}^*\otimes{\C^{p|q}}^*)$ and
$S({\C^d}\otimes\C^{m|n}\oplus{\C^d}^*\otimes{\C^{p|q}}^*)$, these
actions lead to a $gl_d\times gl_{m|n}\times gl_{p|q}$-action on
the latter, where $gl_d$ acts diagonally. It is not immediately
obvious, but nevertheless true \cite{H1}, that
$S({\C^d}\otimes\C^{m|n}\oplus{\C^d}^*\otimes{\C^{p|q}}^*)$ also
admits an action of the larger algebra $gl_d\times gl_{m+p|n+q}$.

For the purpose of describing this action, it is convenient to
introduce a basis for $gl_{m+p|n+q}$ different from that given in
Subsection \ref{Preliminaries}. Set ${\bf I}=\{1, 2, \cdots,
m+p+n+q\}$.  Let $\{v_A |A\in{\bf I}\}$ be a basis of
$\C^{p|q}\oplus\C^{m|n}$ such that $\{v_a |1\le a\le p+q\}$ and
$\{v_{p+q+c} | 1\le c\le m+n\}$ are respectively the standard
bases for $\C^{p|q}$ and $\C^{m|n}$ described in Subsection
\ref{Preliminaries}. Let $E^A_B$, $ A, B\in{\bf I}$,  be the set
of the elementary matrices satisfying $E^A_B v_C=\delta_{B C}
v_A$. These matrices form a homogeneous basis of $gl_{m+p|n+q}$
with the following commutation relations
\begin{eqnarray}
[E^A_B, \ E^C_D]&=&\delta_{B C} E^A_D - (-1)^{{\rm deg}E^A_B{\rm
deg}E^C_ D} \delta_{A D} E^C_B, \label{relation}
\end{eqnarray}
where ${\rm deg}E^A_B$ is the $\Z_2$-degree of $E^A_B$.

Let $B:=\sum_{A\le B; A, B\in{\bf I}} \C E^A_B$, and
$\fh_{m+p|n+q}:=\sum_{A\in{\bf I}} \C E^A_A$. Then $B$ forms a
(non-standard) Borel subalgebra of $gl_{m+p|n+q}$ with the Cartan
subalgebra $\fh_{m+p|n+q}$. Note that $\sum_{a, b=1}^{p+q}\C
E^a_b$ forms a subalgebra isomorphic to $gl_{p|q}$, and $\sum_{u,
v=1}^{m+n}\C E^{p+q+u}_{p+q+v}$ forms a subalgebra isomorphic to
$gl_{m|n}$, and these two subalgebras mutually commute. Together
they form $gl_{p|q}\times gl_{m|n}$, which is a regular subalgebra
of $gl_{m+p|n+q}$ in the sense that its standard Borel subalgebra
$\fb_{p|q}\times\fb_{m|n}$ is contained in $B$, and the
corresponding Cartan subalgebra $\fh_{p|q}\times \fh_{m|n}$ is
identified with $\fh_{m+p|n+q}$. This identification leads to
canonical embeddings of $\fh_{p|q}^*$ and $\fh_{m|n}^*$ in
$\fh_{m+p|n+q}^*$. Choose a basis $\{\widehat\epsilon_A | A\in{\bf
I}\}$ for $\fh_{m+p|n+q}^*$ such that $\widehat\epsilon_A(E_{B
B})=\delta_{A B}$, for all $A, B\in{\bf I}$. An element
$\La=\sum_{A\in{\bf I}} \La_A\widehat\epsilon_A$ of
$\fh_{m+p|n+q}^*$ will also be written as $\La=(\La_1, \La_2,
\cdots, \La_{p+q+m+n})$. Now any pair of elements $\la=(\la_1,
\la_2, \cdots, \la_{p+q})$ $\in$ $\fh_{p|q}^*$ and $\mu=(\mu_1,
\mu_2, \cdots, \mu_{m+n})$  $\in$ $\fh_{m|n}^*$ gives rise to an
element $(\la; \mu) \in \fh_{m+p|n+q}^*$ defined by
\begin{eqnarray}\label{glpqglmnwts}
(\la; \mu)&:=&(\la_1, \la_2,  \cdots, \la_{p+q}, \mu_1, \mu_2,
\cdots, \mu_{m+n}).
\end{eqnarray}

We retain the notations $\C[\x, \etabf]$ and $\C[\y, \zetabf]$
from \secref{finite-duality}, and denote by $\C[\x,\y,\etabf,
\zetabf]$ the polynomial superalgebra $\C[\x,
\etabf]\otimes_\C\C[\y, \zetabf]$. Let $\boldD [\x,
\y,\etabf,\zetabf]$ denote the oscillator superalgebra generated
by the variables ${x^{l}_{i}}$, \ ${\eta^l_{j}}$,\ ${y^{l}_{r}}$,
\ ${\zeta^l_{s}}$, \ and their derivatives
$\frac{\partial}{\partial x^l_{i}}$, \
$\frac{\partial}{\partial\eta^l_{j}}$, \ $\frac{\partial}{\partial
y^l_{r}}$, \ $\frac{\partial}{\partial\zeta^l_{s}}$, where $1\le
i\le m$, \ $1\le j\le n$, \ $1\le r \le p$, \ $1\le s\le q$, and
$1\le l\le d$. Then $\boldD [\x, \y,\etabf,\zetabf]$ naturally
acts on $\C[\x, \y,\etabf,\zetabf]$, thus forms a subalgebra of
${\rm End}(\C[\x, \y,\etabf,\zetabf])$.

The general linear group $GL(d)$ acts on $\boldD[\x,\y,\etabf,
\zetabf]$ by conjugations. The corresponding action of the Lie
algebra $gl_d$ is realized in terms of  the following first order
differential operators($1\le i,j\le d$):
\begin{align}\label{GL}
\Phi(e_{ij})&=\sum_{k=1}^m x_{k}^{i}\frac{\partial}{\partial
x_k^{j}}+\sum_{k=1}^n\eta_k^{i}
\frac{\partial}{\partial\eta_k^{j}}-\sum_{k=1}^p
y_{k}^{j}\frac{\partial}{\partial
y_k^{i}}-\sum_{k=1}^q\zeta_k^{j}\frac{\partial}{\partial\zeta_k^{i}}.
\end{align}
Let $\boldD[\x,\y,\etabf, \zetabf]^{GL(d)}$ denote the
$GL(d)$-invariant subalgebra of $\boldD[\x,\y,\etabf, \zetabf]$.
The $GL(d)$-action is semi-simple.  Thus from the first
fundamental theorem of the invariant theory of the general linear
group (see Chapter 4  of \cite{GW}) we deduce that
$\boldD[\x,\y,\etabf, \zetabf]^{GL(d)}$ is generated by the
following operators :
\begin{align}
\begin{array}{l l}\label{gl{pq}}
\Phi(E^a_b):={\bar\phi}(E_{a b}),\quad\quad& 1\le a, b\le p+q,
\end{array}\end{align}
\begin{align}\begin{array}{l l}\label{gl{mn}}
\Phi(E^{p+q+u}_{p+q+v}):=\phi(E_{u v}), \quad\quad& 1\le u, v\le
m+n,
\end{array}\end{align}
\begin{align}\begin{array}{l l}\label{dd}
\Phi(E^r_{p+q+i}):=\sum_{l=1}^d \frac{\partial}{\partial
y_{r}^{l}}\frac{\partial}{\partial x_{i}^l},\quad\quad &
\Phi(E^r_{m+p+q+j}):=\sum_{l=1}^d\frac{\partial}{\partial
y_r^l}\frac{\partial}{\partial \eta_j^l},\\
\Phi(E^{p+s}_{p+q+i}):=\sum_{l=1}^d \frac{\partial}{\partial
\zeta_{s}^{l}}\frac{\partial}{\partial x_{i}^l},\quad\quad &
\Phi(E^{p+s}_{m+p+q+j}):=\sum_{l=1}^d\frac{\partial}{\partial
\zeta_s^l}\frac{\partial}{\partial \eta_j^l},
\end{array}\end{align}
\begin{align}\begin{array}{l l}\label{xy}
\Phi(E^{p+q+i}_{r}):=-\sum_{l=1}^d x^l_iy_{r}^{l},\quad\quad &
\Phi(E^{p+q+i}_{p+s}):=\sum_{l=1}^dx^l_i\zeta_s^l,\\
\Phi(E^{m+p+q+j}_{r}):=-\sum_{l=1}^d \eta^l_jy_{r}^{l},\quad\quad
& \Phi(E^{m+p+q+j}_{p+s}):=\sum_{l=1}^d\eta^l_j\zeta_s^l,
\end{array}
\end{align}
where $1\le i\le m$, \  $1\le j\le n$, \ $1\le r\le p$, and $1\le
s\le q$. It is an easy exercise to show that the space spanned by
$\Phi(e_{i j})$, \ $1\le i, j\le d$, and $\Phi(E^A_B)$, $A,
B\in{\bf I}$ is a homomorphic image of $gl_d\times gl_{m+p|n+q}$
in $\boldD[\x,\y,\etabf, \zetabf]$. As every Lie superalgebra map
uniquely extends to a homomorphism of its universal enveloping
algebra, we have an associative superalgebra homomorphism
\begin{eqnarray*}\Phi: {\rm U}(gl_d\times
gl_{m+p|n+q})&\rightarrow& \boldD[\x,\y,\etabf,
\zetabf].\end{eqnarray*} Now by identifying
$S({\C^d}\otimes\C^{m|n}\oplus{\C^d}^*\otimes{\C^{p|q}}^*)$ with
the polynomial superalgebra $\C[\x,\y,\etabf, \zetabf]$, we obtain
an action of $gl_d\times gl_{m+p|n+q}$ on
$S({\C^d}\otimes\C^{m|n}+{\C^d}^*\otimes{\C^{p|q}}^*)$. It can be
extracted from \cite{H1} that this action is semi-simple and
multiplicity free. We state this as a theorem for convenience of
reference.
\begin{thm}\cite{H1}\label{howe}
The pair $(gl_d, gl_{m+p|n+q})$ of Lie (super)algebras forms a
dual reductive pair on $\C[\x,\y,\etabf, \zetabf]$.
\end{thm}
\subsection{Unitarity}\label{unitarity}
We first recall some basic facts about $\ast$-superalgebras and
their unitarizable representations.  A $\ast$-superalgebra is an
associative superalgebra $A$ together with an anti-linear
anti-involution $\omega: A\rightarrow A$. Here we should emphasize
that for any $a, b\in A$, we have $\omega(a b)=
\omega(b)\omega(a)$, where no sign factors are involved. A
$\ast$-superalgebra homomorphism $f: (A, \omega) \rightarrow (A',
\omega')$ is a superalgebra homomorphism obeying $f\circ \omega =
\omega'\circ f$. Let $(A, \omega)$ be a $\ast$-superalgebra, and
let $V$ be a $\Z_2$-graded $A$-module. A Hermitian form $\langle\
\cdot\ | \ \cdot \ \rangle$ on $V$ is said to be contravariant if
$\langle a v | v'\rangle = \langle v | \omega(a) v'\rangle$, for
all $a\in A$, $v, v'\in V$. An $A$-module equipped with a positive
definite contravariant Hermitian form is called a unitarizable
$A$-module. It is clear that any unitarizable $A$-module is
completely reducible.

The oscillator superalgebra $\boldD[\x,\y,\etabf, \zetabf]$ admits
the $\ast$-structure $\omega$ defined by
\begin{align*}
\begin{array}{l l l l } {x^{l}_{i}}\mapsto\frac{\partial}{\partial
x^l_{i}}, \quad& \frac{\partial}{\partial x^l_{i}} \mapsto
{x^{l}_{i}},\quad&
{\eta^l_{j}}\mapsto\frac{\partial}{\partial\eta^l_{j}},\quad&
\frac{\partial}{\partial\eta^l_{j}}\mapsto{\eta^l_{j}},\\
{y^{l}_{r}}\mapsto\frac{\partial}{\partial y^l_{r}} \quad&
\frac{\partial}{\partial y^l_{r}}\mapsto{y^{l}_{r}},\quad&
{\zeta^l_{s}}\mapsto\frac{\partial}{\partial\zeta^l_{s}},\quad&
\frac{\partial}{\partial\zeta^l_{s}}\mapsto{\zeta^l_{s}},
\end{array}\end{align*}
where $1\le i\le m$, \ $1\le j\le n$, \ $1\le r \le p$, \ $1\le
s\le q$, and $1\le l\le d$. There exits a unique contravariant
Hermitian form $\langle\ \cdot\ | \ \cdot \ \rangle$  on
$\C[\x,\y,\etabf, \zetabf]$ with $\langle 1 | 1 \rangle=1$. By
using the `particle number' basis relative to which the operators
$x_i^l\frac{\partial}{\partial x_i^l}$, \
${\eta^l_{j}}\frac{\partial}{\partial\eta^l_{j}}$, \ $
{y^{l}_{r}}\frac{\partial}{\partial y^l_{r}}$, \
${\zeta^l_{s}}\frac{\partial}{\partial\zeta^l_{s}}$,  for all $i,
j, r, s, l$ are simultaneously diagonalizable, one can easily show
that the form $\langle\ \cdot\ | \ \cdot \ \rangle$ is positive
definite. The polynomial superalgebra $\C[\x,\y,\etabf, \zetabf]$ with this
inner product (after completion) is the Fock space of $d(m+p)$ bosonic and $d(n+q)$
fermionic quantum oscillators. When $d=1$, we denote it
by $\cF_{m+p|n+q}$. Then it is clear that for arbitrary $d$ we have
$\C[\x,\y,\etabf, \zetabf]$ $\cong$ $\left(\cF_{m+p|n+q}\right)^{\otimes d}$.
What presented in this paragraph is standard material on Fock spaces, which is part of
the basic ingredients of second quantization.

We now consider a $\ast$-structure of ${\rm U}(gl_d\times
gl_{m+p|n+q})$. We shall regard $gl_d\times gl_{m+p|n+q}$ as
embedded in its universal enveloping algebra. Consider the
anti-linear anti-involution $\sigma$ of  ${\rm U}(gl_d\times
gl_{m+p|n+q})$ defined, for all $1\le a, b\le p+q,$\  $p+q+1\le r,
s\le p+q+m+n,$ and $1\le i, j\le d$,  by
\begin{eqnarray*}
E^a_b \mapsto (-1)^{[a]+[b]} E^b_a, &\quad&
E^r_s \mapsto E^s_r,  \\
E^a_s \mapsto -(-1)^{[a]} E^s_a, &\quad& E^r_b \mapsto -(-1)^{[b]} E^b_r, \\
e_{i j}\mapsto e_{j i},
\end{eqnarray*}
where $[a]=\left\{\begin{array}{l l}0, &1\le a\le p,\\
1, & 1\le a-p\le q. \end{array}\right.$  By direct calculations
we can show that this anti-linear map respects the commutation
relations \eqref{relation}, thus indeed defines an anti-linear
anti-involution of  ${\rm U}(gl_d\times gl_{m+p|n+q})$. Now
$\sigma$ gives rise to a $\ast$-structure for ${\rm U}(gl_d\times
gl_{m+p|n+q})$.
\begin{thm}(1). The map $\Phi$ is a $\ast$-superalgebra
homomorphism from $({\rm U}(gl_d\times gl_{m+p|n+q}), \sigma)$ to
the oscillator superalgebra $(\boldD[\x,\y,\etabf, \zetabf], \omega)$.\\
(2). $\C[\x, \y,\etabf,\zetabf]$ is a unitarizable $({\rm
U}(gl_d\times gl_{m+p|n+q}), \sigma)$-module with respect to the
Hermitian form $\langle \, \cdot\, , \, \cdot \, \rangle$.
\end{thm}
\begin{proof} Using  equations \eqref{GL}-\eqref{xy},
we can show by direct calculations that for all $X\in gl_d\times
gl_{m+p|n+q}$, we have $\Phi\sigma(X) = \omega\Phi(X).$ This
proves part (1). Part (2) immediately follows from part (1).
\end{proof}
We also have the following result.
\begin{lem}\label{hwtype} All the irreducible $gl_d\times
gl_{m+p|n+q}$-submodules of $\C[\x, \y,\etabf,\zetabf]$ are of
highest weight type with respect to the Borel subalgebra
$\fb_d\times B$.
\end{lem}
\begin{proof}
Let $H$ be the harmonic subspace of $\C[\x, \y,\etabf,\zetabf]$,
i.e., the subspace consisting of such polynomials that are annihilated by
all the elements $\Phi(E^A_B)$, $A\le p+q$, $B>p+q$, of
\eqref{dd}. Then $H$ forms a module of the subalgebra $gl_d\times
gl_{p|q}\times gl_{m|n}$ spanned by elements of (\ref{GL}),
\eqref{gl{pq}} and \eqref{gl{mn}}. It follows from the $(gl_d,
gl_{m|n})$-duality on $\C[\x, \etabf]$ described in
\thmref{hw-glmn} and the $(gl_d, gl_{p|q})$-duality on $\C[\y,
\zetabf]$ described in \thmref{hw*} that $H$ decomposes into a
direct sum of finite dimensional irreducible $gl_d\times
gl_{p|q}\times gl_{m|n}$-modules.

Let $W$ be an irreducible $gl_d\times gl_{m+p|n+q}$-submodule of
$\C[\x, \y,\etabf,\zetabf]$. Let $H_W=W\cap H$. Then $H_W\ne 0$,
as the lowest order polynomials of $W$ are all contained in $H_W$.
Now $H_W$ forms a module of the subalgebra $gl_d\times
gl_{p|q}\times gl_{m|n}$, which in fact is irreducible.  To see
this, we note that if $H_W$ were reducible with respect to
$gl_d\times gl_{p|q}\times gl_{m|n}$, then due to complete
reducibility $H_W$ would contain more than one linearly
independent $gl_d\times gl_{p|q}\times gl_{m|n}$-highest weight
vectors. However, since they lie in $H_W$, they would also be
$gl_d\times gl_{m+p|n+q}$-highest weight vectors with respect to
$\fb_d\times B$, thus contradicting the irreducibility of $W$.
Thus $W$ is generated by a $gl_d\times gl_{m+p|n+q}$-highest
weight vectors with respect to $\fb_d\times B$, as claimed.
\end{proof}

\begin{rem}
From  \thmref{howe} and the proof of the lemma we can deduce that
the action of $gl_d\times gl_{p|q}\times gl_{m|n}$ on the harmonic
subspace $H$ of $\C[\x, \y,\etabf,\zetabf]$ is semi-simple and
multiplicity free.
\end{rem}

Let $\fg^\R$ be the real superspace
spanned by $\{X\in (gl_{m+p|n+q})_{\bar{0}}
| \sigma(X)=-X\}\cup \sqrt{i}\{X\in (gl_{m+p|n+q})_{\bar{1}} |
\sigma(X)=-X\}$. Then $\fg^\R$ is a real form of $gl_{m+p|n+q}$, that is, $\fg^\R$
forms a real Lie superalgebra with the complexification being
$gl_{m+p|n+q}$ itself. The usual notation for this real form is
$\fu(m, p|n, q)$. Note that the maximal even subalgebra of
$\fg^\R$ is $\fu(m, p)\times\fu(n, q)$. Thus every nontrivial
unitarizable $\fu(m, p|n, q)$-module  must be infinite
dimensional.

\subsection{Comments on unitarizable modules}\label{comments}
At this point, we should relate to results in the literature. Note
that the restrictions of $\sigma$ to the subalgebras $gl_{m|n}$
and $gl_{p|q}$ act differently on the odd subspaces. They
respectively give rise to two different real forms $u_+(m|n)$ and
$u_-(p|q)$ of the subalgebras. Now $\fg^\R$ contains the
subalgebra $\fu_-(p|q)\times\fu_+(m|n)$, which one would like to
regard as the `maximal compact subalgebra'. Unfortunately finite
dimensional representations of $\fu_+(m|n)$ and $\fu_-(p|q)$ are
not necessarily unitarizable. In fact it has long been known
\cite{GZ} that the only finite dimensional unitarizable
irreducible representations of $u_+(m|n)$ (resp. $u_-(p|q)$) are
the tensor products of the irreducible representations appearing
in \thmref{hw-glmn} (resp. \thmref{hw*}) with some $1$-dimensional
representations (upon restricting modules of the general linear
superalgebra to modules of its real form), which constitute a
small class of the finite dimensional irreducible representations.
Thus the situation is very different from the case of the compact
real Lie algebra $\fu(k)$.

The intersection of $\fu(m, p|n, q)$ with $sl_{m+p|n+p}$ gives
rise to the real Lie superalgebra $su(m, p|n, q)$. It was shown in
\cite{J} that $su(m, p|n, q)$ admits no unitarizable highest or
lowest weight representations with respect to the standard
Borel subalgebra if all the integers $m$,  $n$, $p$ and $q$
are non-zero. Since \cite{J} was only concerned with irreducible unitarizable
highest weight modules of simple basic classical Lie superalgebras
with respect to their standard Borel subalgebras \cite{K},
the irreducible $gl_{m+p|n+q}$-modules appearing
in the decomposition of $\C[\x, \y,\etabf,\zetabf]$ were ignored.
In fact with respect to the standard Borel subalgebra of $gl_{m+p|n+q}$,
the unitarizable irreducible representations studied in this paper are neither
highest weight nor lowest weight type unless some of the integers
$m$,  $n$, $p$ and $q$ are zero.

A final comment is that when both of the integers $n$ and $q$ are
zero, the general linear superalgebra $gl_{m+p|n+q}$ reduces to
the ordinary Lie algebra $gl_{m+p}$, and $\C[\x,
\y,\etabf,\zetabf]$ to the ordinary polynomial algebra $\C[\x,
\y]$ in the two sets of variables $\x$ and $\y$. Then the
irreducible $gl_{m+p}$-modules appearing in $\C[\x, \y]$ are the
unitarizable irreducible $u(m, p)$-modules studied by Kashiwara
and Vergne in \cite{KV}. It is known \cite{DES, KV} that the
unitarizable irreducible $u(m, p)$-module at every reduction point
\cite{EHW} is a submodule in $\C[\x, \y]$ for some $d$. However,
it is not known whether this is still true in the super case.

\subsection{The $(gl_d,gl_{m+p|n+q})$-duality on
$\C[\x, \y,\etabf,\zetabf]$}\label{glglduality**}
Each generalized partition $\la=(\la_1,\cdots,\la_d)$ of length $d$ can be
uniquely expressed as $\la=\la^+ + \la^-$, with
\begin{eqnarray*}\la^+&:=&(\hbox{max}\{\la_1,0\},\cdots,\hbox{max}\{\la_d,0\}),\\
\la^-&:=&(\hbox{min}\{\la_1,0\},\cdots,\hbox{min}\{\la_d,0\}).
\end{eqnarray*}
Note that $\la^+$ is a partition of length $d$ , while $\la^-$ is a generalized
partition of non-positive integers with length $d$, such that $(\lambda^-)^*$
is a partition.  Furthermore,
\begin{eqnarray} \label{depthlimit} \mbox{depth of $\la^+$}
\, + \ \mbox{depth of $(\la^-)^*$}&\le & d,
\end{eqnarray}
where the depth of a partition is the number of positive integers in
it.

The generalized partition $\la=(\la_1,\cdots,\la_d)$ satisfies the
conditions $\la_{m+1}\le n$ and $\la_{d-p}\ge -q$ if and only if
$\la^+_{m+1}\le n$ and $(\la^-)^*_{p+1}\le q$. Corresponding to
each such generalized partition $\la$, we define
\begin{equation}\label{Deltala**}
\Box_\la:= \Delta^*_{\la^-}\, \Delta_{\la^+}.
\end{equation}
\begin{lem} \label{Box}
If the generalized partition $\la$ satisfies the conditions
$\la_{m+1}\le n$ and $\la_{d-p}\ge -q$, then $\Box_\la$ is a
non-zero $gl_d \times gl_{m+p|n+q}$ highest weight vector with
respect to the Borel subalgebra $\fb_d\times B$. The $gl_d$-weight
of $\Box_\la$ is $\la$, and the $gl_{m+p|n+q}$-weight is given by
\begin{eqnarray}\label{Lala}
\La(\la)&:=&(-d {\bf 1}+ \widehat{\la^-};\  \widetilde{\la^+}),
\end{eqnarray}
where the expression on the right hand side is as explained by
\eqref{glpqglmnwts}.
\end{lem}
\begin{proof}
By construction, $\Box_\la$ is a highest weight vector of the
subalgebra $gl_d\times gl_{p|q}\times gl_{m|n}$ with respect to
the standard Borel subalgebra $\fb_d\times
\fb_{p|q}\times\fb_{m|n}$. Therefore, we only need to show that
$\Phi(E_{p+q, p+q+1})= \sum_{k=1}^d
\frac{\partial}{\partial\zeta^k_q}\frac{\partial}{\partial x^k_1}$
annihilates $\Box_\la$ in order to prove that $\Box_\la$ is a
$gl_d \times gl_{m+p|n+q}$ highest weight vector with respect to
the Borel subalgebra $\fb_d\times B$. If $\Phi(E_{p+q,
p+q+1})\Box_\la\ne 0$, then there must exist at least one integer
$i\in\{1, 2, \cdots, d\}$ such that $x_1^i \zeta_q^i$ appears in
$\Box_\la$. Let $ht(\la^+)$ denote the depth of $\la^+$, and
$ht((\la^-)^*)$ denote the depth of $(\la^-)^*$. By examining
its explicit form,  we can see that $\Box_\la$ does not involve
any of the variables $x_1^k$, $d\ge k>ht(\la^+)$,
and $\zeta_q^l$, $1\le l< d+1-ht((\la^-)^*)$.
Therefore in order for $x_1^i \zeta_q^i$ to
appear in $\Box_\la$, the integer $i$ must satisfy
$d+1-ht((\la^-)^*)\le i\le ht(\la^+)$. But this is impossible
since \eqref{depthlimit} requires $ht(\la^+)+ht((\la^-)^*)\le d$.
Therefore, $\Phi(E_{p+q, p+q+1})\Box_\la=0,$ and thus $\Box_\la$
is a $gl_d \times gl_{m+p|n+q}$ highest weight vector with respect
to the Borel subalgebra $\fb_d\times B$.

The $gl_d$-weight of $\Box_\la$ is obviously $\la$.  From
\lemref{glpmn-duality} and \lemref{glpq-hwv}, we easily see that
the $gl_{m+p|n+q}$-weight of $\Box_\la$ is indeed $(-d {\bf 1}+
\widehat{\la^-};\  \widetilde{\la^+})$.
\end{proof}

We shall denote by $W^{\La}_{m+p|n+q}$ the irreducible highest
weight $gl_{m+p|n+q}$-module with highest weight $\La$ relative to
the non-standard Borel subalgebra $B$.
\begin{thm}\label{hw+*}
Under the $gl_d\times gl_{m+p|n+q}$-action $\C[\x,\y,\etabf,
\zetabf]$ decomposes into
\begin{equation}
\C[\x,\y,\etabf, \zetabf]\cong\sum_{\la}V^\la_d\otimes
W^{\La(\la)}_{m+p|n+q},\label{glgl**}
\end{equation}
where $\la$ is summed over all generalized partitions of length
$d$ subject to $\la_{m+1}\le n$ and $\la_{d-p}\ge -q$.
Furthermore, $\C\Box_\lambda$ is the space of highest weight
vectors of the irreducible module $V^\la_d\otimes
W^{\La(\la)}_{m+p|n+q}$.
\end{thm}
\begin{proof} Note that every irreducible $gl_d$-submodule of
$\C[\x,\y,\etabf, \zetabf]$ is finite dimensional. Thus it follows
from \thmref{howe} and \lemref{hwtype} that the decomposition of
$\C[\x,\y,\etabf, \zetabf]$ under $gl_d\times gl_{m+p|n+q}$ has to
be of the form (\ref{glgl**}), with the sum in $\la$ ranging over
some subset of generalized partitions of length $d$. (Here if
$\la$ happens to be a generalized partition not satisfying
$\la_{m+1}\le n$ and $\la_{d-p}\ge -q$, the expression
$\Lambda(\la)$ stands for the highest weight for $gl_{m|p|n|q}$
corresponding to $\la$ under this Howe duality.) In view of
\lemref{Box}, we only need to show that every generalized
partition $\la$ belonging to this subset must satisfy the
conditions $\la_{m+1}\le n$ and $\la_{d-p}\ge -q$,  in order to
prove the theorem.

Let $\mu$ be a generalized partition of length $d$.
Assume that either one or both of the conditions $\mu_{m+1}\le
n$ and $\mu_{d-p}\ge -q$ are violated. We choose a pair of
positive integers $p'$ and $n'$ with $p'\ge p$ and $n'\ge n$ such
that $\mu_{m+1}\le n'$ and $\mu_{d-p'}\ge -q$. Such an $n^\prime$
is trivial to come by, and such a $p^\prime$ also exists since by Remark \ref{pged}
the condition $\mu_{d-p'}\ge -q$ is always satisfied if $p^\prime\ge
d$. We let $\C[\x,\bar\y,\bar\etabf, \zetabf]$ denote the
polynomial superalgebra generated by
\begin{equation*}
x_i^l:=e^l\otimes e_i,\quad \bar y_r^l:={e^l}^*\otimes e_r^*,\quad
\bar\eta_j^l:=e^l\otimes f_j,\quad \zeta_s^l:={e^l}^*\otimes
f_s^*,
\end{equation*}
for $1\le l\le d$, $1\le i\le m$, $1\le j\le n'$, $1\le r\le p'$
and $1\le s\le q$.  Then $\C[\x,\y,\etabf, \zetabf]$ becomes a
subspace of $\C[\x,\bar\y,\bar\etabf, \zetabf]$ upon identifying
\begin{equation}
y_r^l:={\bar y}_{p^\prime-p+r}^l,\quad \zeta_s^l:=\zeta_s^l, \quad
x_i^l:=x_i^l,\quad \eta_j^l:=\bar\eta_j^l, \label{identity}
\end{equation}
for $1\le l\le d$, $1\le i\le m$, $1\le j\le n$, $1\le r\le p$ and
$1\le s\le q$. We denote this inclusion by $\iota:
\C[\x,\y,\etabf, \zetabf]\rightarrow \C[\x,\bar\y,\bar\etabf,
\zetabf]$. There is also the surjection $\pi:
\C[\x,\bar\y,\bar\etabf, \zetabf]\rightarrow \C[\x,\y,\etabf,
\zetabf]$ defined by setting ${\bar y}_r^l=0$, $1\le r\le p^\prime
-p$, and ${\bar\eta}_{s}^l=0$, $n<s\le n^\prime$, for all $l$,
then making the identification (\ref{identity}). Obviously,
$\pi\iota$ is the identity map on $\C[\x,\y,\etabf, \zetabf].$

Now we turn to the $gl_d\times gl_{m+p'|n'+q}$-action on
$\C[\x,\bar\y,\bar\etabf, \zetabf]$. Upon choosing the basis for
$\C^{p'|q}\oplus\C^{m|n'}$ that is the union of the standard bases
of $\C^{p'|q}$ and $\C^{m|n'}$, the general linear superalgebra
$gl_{m+p'|n'+q}$ becomes the Lie superalgebra of
$(m+p'+n'+q)\times(m+p'+n'+q)$-matrices. Consider the subalgebra
$\fl$ of $gl_{m+p'|n'+q}$ consisting of matrices of the form
$\left(\begin{array}{c c c} D & 0 & 0\\
0 & X & 0
\\ 0 & 0 & D'\end{array}\right)$, where
$X{\in} gl_{m+p|n+q}$, and $D$ and $D'$ are diagonal matrices of
sizes $(p^\prime -p)\times(p^\prime -p)$ and $(n'-n)\times(n'-n)$,
respectively. Obviously $\fl$ contains the $gl_{m+p|n+q}$
subalgebra
$\left\{\left.\left(\begin{array}{c c c} 0 & 0 & 0\\
0 & X & 0
\\ 0 & 0 & 0\end{array}\right)\right|
X{\in}gl_{m+p|n+q}\right\}.$  Let $\fp=\fn + \fl$ be a parabolic
subalgebra of $gl_{m+p'|n'+q}$ with the Levi factor $\fl$ and
nilpotent radical $\fn$. We assume that $\fp$ contains all the
upper triangular matrices.  Then there exists a nilpotent
subalgebra $\bar\fn$ consisting of strictly lower triangular
matrices such that $gl_{m+p'|n'+q}=\fp + \bar\fn$. By examining
equations \eqref{GL}-\eqref{xy} we can see that $\iota$ is a
$gl_d\times gl_{m+p|n+q}$-module map. Let $V$  be any $gl_d\times
gl_{m+p|n+q}$-submodule of $\C[\x, \y, \etabf, \zetabf]$. Then
$\iota(V)$ is in fact a $gl_d\times\fp$-module with $\fn$ acting
by zero. Thus $W=\Phi({\rm U}(\bar\fn))\iota(V)$ forms a
$gl_d\times gl_{m+p'|n'+q}$-submodule of $\C[\x,\bar\y,\bar\etabf,
\zetabf]$. Note that $W$ is irreducible if $V$ is irreducible with
respect to $gl_d\times gl_{m+p|n+q}$. Again by examining equations
\eqref{GL}-\eqref{xy} we can see that $\pi$ is a $gl_d\times
gl_{m+p|n+q}$-module map from the restriction of
$\C[\x,\bar\y,\bar\etabf, \zetabf]$ to $\C[\x, \y, \etabf,
\zetabf]$, and satisfies $\pi(W)=V$. This in particular implies
that $\pi\iota$ is the identity $gl_d\times gl_{m+p|n+q}$-module
map on $\C[\x,\y,\etabf, \zetabf].$

Let $v_\mu\in \C[\x, \y, \etabf, \zetabf]$ be any $gl_d\times
gl_{m+p|n+q}$ highest weight vector with the $gl_d$-weight
$\mu$ (that violates one or both of the conditions $\mu_{m+1}\le n$
and $\mu_{d-p}\ge -q$). Then by the above discussion,
$\iota(v_\mu)$ is a $gl_d\times\fl$ highest weight vector in $
\C[\x,\bar y,\bar\etabf, \zetabf]$, which has the same $gl_d$ weight,
and is also annihilated by $\fn$. Therefore, $\iota(v_\mu)$ is a
$gl_d \times gl_{m+p'|n'+q}$ highest weight vector in $ \C[\x,\bar
y,\bar\etabf, \zetabf]$. By Theorem \ref{howe} and Lemma \ref{Box}
(with $p$ replaced by $p'$ and $n$ by $n'$), there exists a unique
non-zero $\Box_\mu\in \C[\x,\bar\y,\bar\etabf, \zetabf]$ such that
$\iota(v_\mu)=c \Box_\mu$ for some complex number $c$.

We claim that every monomial in the polynomial $\Box_\mu$ contains
at least one of the variables $\bar y_r^l$, $\bar\eta_s^l$, where
$r=1,\cdots,p'-p$, \ $s=n+1,\cdots,n'$ and $l=1,\cdots,d$. This
can be seen from the explicit form \eqref{Deltala**} of
$\Box_\la$. Consider first the case with $t:=\mu_{m+1}>n$. Then
$(\mu^+)'_1>m$, and $\Delta_{\mu^+}$ has the factor $\Delta_{t,
(\mu^+)'_t}$. From \eqref{eqdet} we see that $\Delta_{t,
(\mu^+)'_t}$ is the determinant of a matrix with rows
$(\bar\eta^1_t, \bar\eta^2_t, \cdots, \bar\eta^{(\mu^+)'_t}_t)$
with Grassmann number entries. Now consider the case with
$\mu_{d-p}<-q$. Let $\gamma=(\mu^-)^*$, then $\gamma_{p+1}>q$.
Thus we must also have $\gamma_1>q$, and
$u:=\gamma^\prime_{q+1}\ge p+1$. Now $\Delta^*_{\mu^-}$ has the
factor $\Delta^*_u$.  From \eqref{eqdet*} we see that $\Delta^*_u$
is the determinant of a matrix with a column $(\bar y^d_{p'+1-u},
\bar y^{d-1}_{p'+1-u}, \cdots, \bar y^{d+1-u}_{p'+1-u})$. Note
that $p'+1-u\le p'-p$.

Now it is obvious that $\pi(\Box_\mu)=0$, which in turn implies $v_\mu=0$.
Therefore, the decomposition of $\C[\x, \y, \etabf, \zetabf]$ can
not contain $V^\mu_d\otimes W^{\La(\mu)}_{m+p|n+q}$ as an
irreducible submodule if $\mu$ violates any of the conditions
$\mu_{m+1}\le n$ and $\mu_{d-p}\ge -q$.
\end{proof}

 By using \thmref{hw-glmn}, \thmref{hw*} and the decomposition $\C[\x, \y, \etabf, \zetabf]=\C[\x,
\etabf]\otimes_\C\C[\y, \zetabf]$, we have
\begin{eqnarray*}
\C[\x, \y, \etabf, \zetabf]&\cong&\sum_{\la, \mu} V_d^{\la}\otimes
V_d^{\mu}\otimes V^{\tilde\la}_{m|n}\otimes V^{-d {\bf 1}+
\hat\mu}_{p|q},
\end{eqnarray*}
where the summation in $\la$ is over all the partitions of length
$d$ satisfying $\la_{m+1}\le n$, and the summation in $\mu$ is
over all the generalized partitions of non-positive integers of
length $d$ satisfying $\mu_{d-p}\ge -q$.

The decomposition of the tensor product of any two finite
dimensional irreducible $gl_d$-modules is described by the
Littlewood-Richardson theory. For any generalized partitions $\la$
and $\mu$ of length $d$,
\begin{eqnarray}
V_d^{\la}\otimes V_d^{\mu}&\cong&\sum_{\nu} C_{\la \mu}^\nu
V_d^{\nu}, \label{LR}
\end{eqnarray}
where the non-negative integers $C_{\la \mu}^\nu$ are the
so-called Littlewood-Richardson coefficients, which give the
respective multiplicities of the irreducible $gl_d$-modules
$V_d^{\nu}$ appearing in the tensor product.

Therefore, $\C[\x, \y, \etabf, \zetabf]$ decomposes into
\begin{eqnarray}\label{restriction}
\C[\x, \y, \etabf, \zetabf]&\cong&\sum_\nu V_d^{\nu}\otimes
\sum_{\la, \mu} C_{\la \mu}^\nu V^{\tilde\la}_{m|n}\otimes V^{-d
{\bf 1}+ \hat\mu}_{p|q},
\end{eqnarray}
where the summation in $\nu$ is over all the generalized
partitions of length $d$. Combining Theorem \ref{hw+*} and
\eqref{restriction}, the Littlewood-Richardson coefficients
$C_{\la \mu}^\nu$ appearing in \eqref{restriction} may be non-zero
only when $\nu$ satisfies the conditions $\nu_{m+1}\le n$ and
$\nu_{d-p}\ge -q$. Now if we put $p=d-m-1$, $n=\la_{m+1}$ and
$q=-\mu_{d-p}$, then the conditions become $\nu_{m+1}\le
\la_{m+1}$ and $\nu_{m+1}\ge \mu_{m+1}$. Letting $m$ run from $0$
to $d-1$, we have the following corollary.

\begin{cor}\label{coefficient1}
 Assume that $\la$ is a partition of length $d$ and
$\mu$ is a generalized partition of non-positive integers of
length $d$. If $\nu$ is a generalized partition of length $d$
satisfies $\nu_{m}> \la_{m}$ or $\nu_{m}< \mu_{m}$ for some $m
\in\{1,2,\dots,d\}$, then $C_{\la \mu}^\nu=0$.
\end{cor}

\corref{coefficient1} translated to ordinary partitions implies
the following.

\begin{cor}\label{coefficient2} Let $\la$ and $\mu$ be two partitions of length $d$. If $\nu$ is
a partition of length $d$ satisfies $\nu_{m}> \rm{
min}\{\mu_m+\la_1, \la_m+\mu_1\}$ for some $m \in\{1,2,\dots,d\}$,
then the Littlewood-Richardson coefficient $C_{\la \mu}^\nu=0$.
\end{cor}

\begin{proof} Let $\la=(\la_1,\la_2,\dots,\la_d)$ and
$\mu=(\mu_1,\mu_2,\dots\mu_d)$ be two partitions. Then
$\mu-\mu_1{\bf 1}:=(\mu_1-\mu_1,\mu_2-\mu_1,\dots\mu_d-\mu_1)$ is
a generalized partition. By \corref{coefficient1}, we have $C_{\la
\mu}^\nu=C_{\la, \mu-\mu_1{\bf 1}}^{\nu-\mu_1{\bf 1}}=0$ if
$\nu_{m}-\mu_1> \la_m$ for some $m \in\{1,2,\dots,d\}$. Therefore,
the corollary follows from the symmetry property of the
Littlewood-Richardson coefficients.
\end{proof}

\begin{rem} A alternative method to prove \thmref{hw+*} is the
following. One can first prove \corref{coefficient2}, using, for
example, the celebrated combinatorial algorithm known as the
Littlewood-Richardson rule (see e.g.~\cite{M}). Using
\corref{coefficient2} it can then be derived that in the tensor
product decomposition of $V^\la_d\otimes V^\mu_d$, with $\la$ a
partition of length $d$ satisfying $\la_{m+1}\le n$, and $\mu$ a
generalized partition of non-positive integers of length $d$
satisfying $\mu_{d-p}\ge -q$, only $gl_d$-modules associated to
generalized partitions $\nu$ with $\nu_{m+1}\le n$ and
$\nu_{d-p}\ge -q$ can occur.  Using this fact together with
\lemref{Box} it is then not difficult to prove \thmref{hw+*}.
\end{rem}

\section{Branching rules of unitarizable irreducible
$gl_{m+p|n+q}$-modules}\label{branching}

As an easy application of \thmref{hw+*}, we derive the
$gl_{m+p|n+q}\rightarrow gl_{p|q}\times gl_{m|n}$ branching rule
for the infinite dimensional unitarizable irreducible
$gl_{m+p|n+q}$-representations arising from the decomposition of
tensor powers of the Fock space of $m+p$ bosonic and $n+q$
fermionic quantum oscillators. Results of this section will be
important for developing a character formula for these
unitarizable irreducible $gl_{m+p|n+q}$-modules.

Let us denote by $\fk^\C$ the subalgebra $gl_{p|q}\times gl_{m|n}$
of $gl_{m+p|n+q}$. Recall the decomposition of $\C[\x, \y, \etabf,
\zetabf]$ as a $gl_d\times\fk^\C$-module \eqnref{restriction}.
 Denote by
$\left.W_{m+p|n+q}^{\La(\nu)}\right|_{\fk^C}$ the restriction of
$W_{m+p|n+q}^{\La(\nu)}$ to a $\fk^\C$-module. Let us now consider
Theorem \ref{hw+*} by restricting both sides of equation
\eqref{glgl**} to $gl_d\times \fk^\C$-modules. Using
\eqref{restriction} we obtain
\begin{eqnarray*}
\sum_\nu V_d^{\nu}\otimes
\left.W_{m+p|n+q}^{\La(\nu)}\right|_{\fk^C}&\cong&\sum_\nu
V_d^{\nu}\otimes \sum_{\la, \mu} C_{\la \mu}^\nu
V^{\tilde\la}_{m|n}\otimes V^{-d {\bf 1}+ \hat\mu}_{p|q},
\end{eqnarray*}
where the summation in $\nu$ on the left hand side is over the
generalized partitions of length $d$ satisfying the conditions
$\nu_{m+1}\le n$ and $\nu_{d-p}\ge -q$. The above equation
immediately leads to the following $gl_{m+p|n+q}\rightarrow
gl_{p|q}\times gl_{m|n}$ branching rule:

\begin{thm} Let $\nu$ be a generalized partition of length $d$
subject to the conditions $\nu_{m+1}\le n$ and $\nu_{d-p}\ge -q$.
We have
\begin{eqnarray}\label{branchrule}
\left.W_{m+p|n+q}^{\La(\nu)}\right|_{gl_{p|q}\times
gl_{m|n}}&\cong&\sum_{\la, \mu} C_{\la \mu}^\nu
V^{\tilde\la}_{m|n}\otimes V^{-d {\bf 1}+ \hat\mu}_{p|q},
\end{eqnarray}
where the summation in $\la$ is over all the partitions of length
$d$ satisfying $\la_{m+1}\le n$, and the summation in $\mu$ is
over all the generalized partitions of non-positive integers with
length $d$ satisfying $\mu_{d-p}\ge -q$.
\end{thm}

Recall that $\La(\nu)$ is defined by \eqref{Lala}.

\section{Character formula for unitarizable irreducible
$gl_{m+p|n+q}$-modules} \label{character-formula} In this section,
we shall develop a character formula for the infinite dimensional
unitarizable irreducible $gl_{m+p|n+q}$-modules appearing in the
decomposition of $\C[\x, \y, \etabf, \zetabf]$. Let us first
present some background material on Schur functions and the
so-called hook Schur functions of Berele-Regev \cite{BR}. A
comprehensive reference on Schur functions is \cite{M}.

\subsection{Hook Schur function}   Let $\x=\{x_1,x_2,\cdots ,x_m\}$
be a set of $m$ variables. To a partition $\la$ of non-negative
integers we may associate the Schur function $s_\la(x_1,x_2,\cdots
, x_m)$. We will write $s_\la(\x)$ for $s_\la(x_1,x_2,\cdots
,x_m)$.  For a partition $\mu\subset\la$ we let $s_{\la/\mu}(\x)$
denote the corresponding skew Schur function. Denote by $\mu'$ the
conjugate partition of a partition $\mu$.  The {\em hook Schur
function} \cite{BR} corresponding to a partition $\la$ is defined
by
\begin{equation}\label{hookschur1}
HS_{\la}(\x;\y):=\sum_{\mu\subset\la}s_\mu(\x)s_{\la'/\mu'}(\y),
\end{equation}
where as usual $\y=\{y_1,y_2,\cdots ,y_n\}$.

Let $\la$ be a partition and $\mu\subseteq\la$.  We fill the boxes
in $\mu$ with entries from the linearly ordered set
$\{x_1<x_2<\cdots <x_m\}$ so that the resulting tableau is
semi-standard. Recall that this means that the rows are
non-decreasing, while the columns are strictly increasing.  Next
we fill the skew partition $\la/\mu$ with entries from the
linearly ordered set $\{y_1<y_2<\cdots <y_n\}$ so that it is
conjugate semi-standard, which means that the rows are strictly
increasing, while its columns are non-decreasing.  We will refer
to such a tableau as an {\em $(m|n)$-semi-standard tableau}
(cf.~\cite{BR}).  To each such tableau $T$ we may associate a
polynomial $(\x\y)^T$, which is obtained by taking the products of
all the entries in $T$. Then we have \cite{BR}
\begin{equation}\label{hookschur}
HS_{\la}(\x;\y)=\sum_{T}(\x\y)^T,
\end{equation}
where the summation is over all $(m|n)$-semi-standard tableaux of
shape $\la$.

We recall the following combinatorial identity involving hook
Schur functions that is of crucial importance in the sequel. Note
that $HS_{\la}(\x;\y)\not=0$ iff $\la_{m+1}\le n$.

\begin{prop}\cite{CL}\label{combid}
Let $\x=\{x_1,x_2,\cdots ,x_m\}$, $\etabf=\{\eta_1,\eta_2,\cdots
,\eta_n\}$ be two sets of variables and ${\bf z}=\{z_1,z_2,\cdots,z_d\}$
be $d$ variables. Then
\begin{equation*}
\prod_{i,j,k}(1-x_iz_k)^{-1}(1+\eta_jz_k)=\sum_{\la}HS_{\la}(\x;\etabf)s_{\la}(\bf
z),
\end{equation*}
where $1\le i\le m$, $1\le j\le n$, $1\le k\le d$ and $\la$ is
summed over all partitions with length $d$ subject to
$\la_{m+1}\le n$.
\end{prop}

Replacing $x_i$, $\eta_j$ and $z_k$ in \propref{combid} by
$y_i^{-1}$, $\zeta_j^{-1}$ and $z_k^{-1}$, we obtain the following
result.

\begin{prop}\cite{CL}\label{combid*}
Let $\y=\{y_1,y_2,\cdots ,y_p\}$,
$\zetabf=\{\zeta_1,\zeta_2,\cdots ,\zeta_q\}$,
and ${\bf z}=\{z_1,z_2,\cdots,z_d\}$. Set
$\y^{-1}=\{y^{-1}_1,y^{-1}_2,\cdots ,y^{-1}_p\}$,
$\zetabf^{-1}=\{\zeta^{-1}_1,\zeta^{-1}_2,\cdots ,\zeta^{-1}_q\}$,
and ${\bf z}^{-1}=\{z^{-1}_1,z^{-1}_2,\cdots,z^{-1}_d\}$.  Then
\begin{equation*}
\prod_{i,j,k}(1-y_i^{-1}z_k^{-1})^{-1}(1+\zeta_j^{-1}z_k^{-1})
=\sum_{\la}HS_{\la}(\y^{-1};\zetabf^{-1})s_{\la}(\bf z^{-1}),
\end{equation*}
where $1\le k\le d$, $1\le i\le p$ and $1\le j\le q$ and $\la$ is
summed over all partitions with length $d$ subject to
$\la_{p+1}\le q$.
\end{prop}

We recall the following lemma which plays a crucial role in
developing a character formula for unitarizable irreducible
$gl_{m+p|n+q}$-modules by using Howe duality. Denote by
$ch(V_d^\la)$ the formal character of the irreducible
$gl_d$-module $V_d^\la$.

\begin{lem}\cite{CL}\label{independence}
Let $q$ be an indeterminate and suppose that $\sum_\la\phi(q){\rm
ch}V^\la_d=0$, where $\phi_\la(q)$ are power series in $q$ and
$\la$ above is summed over all generalized partitions of length
$d$. Then $\phi_\la(q)=0$, for all $\la$.
\end{lem}

\subsection{Character formula for finite dimensional modules}
Recall from Subsection \ref{Preliminaries} that
$\tilde\epsilon_1,\cdots,\tilde\epsilon_{d}$, are the weights of
the natural $gl_d$-module $\C^d$, and $\epsilon_1, \cdots,
\epsilon_m$, $\delta_1,\cdots,\delta_n$ are  the weights of the
natural $gl_{m|n}$-module $\C^{m|n}$. Let $e$ be a formal
indeterminate. For $k=1,\cdots,d$,\quad $i=1,\cdots, m$ and
$j=1,\cdots, n$, we set
\begin{equation}
\bar z_k=e^{\tilde\epsilon_k}, \quad \bar x_i=e^{\epsilon_i},\quad
\bar\eta_j=e^{\delta_j},
\end{equation}
and let $\bar\x=\{\bar x_1,\bar x_2,\cdots ,\bar x_m\}$, $\bar
\etabf=\{\bar \eta_1,\bar \eta_2,\cdots ,\bar \eta_n\}$ and $\bar
{\bf z}=\{\bar z_1,\bar z_2,\cdots,\bar z_d\}$.

Consider $\C[\x,\etabf]$ as a $gl_d\times gl_{m|n}$-module. Its
formal character ${\rm ch}(\C[\x,\etabf])$ with respect to the
Cartan subalgebra $\sum_{i=1}^{m+n} \C E_{i i}
\oplus\sum_{k=1}^d\C e_{k k}$ can be easily computed by using
equations \eqref{glpq} and \eqref{glmn1} to give
\begin{equation}
{\rm ch}(\C[\x,\etabf]) =\prod_{i,j,k}(1-\bar x_i \bar
z_k)^{-1}(1+\bar \eta_j\bar z_k),
\end{equation}
where $1\le i\le m$, $1\le j\le n$, $1\le k\le d$.  Thus by
\propref{combid}
\begin{equation}\label{character1}
{\rm ch}(\C[\x,\etabf]) =\sum_{\la}HS_{\la}(\bar \x;\bar
\etabf)s_{\la}(\bar{\bf z}),
\end{equation}
where $\la$ is summed over all partitions with length $d$ subject
to $\la_{m+1}\le n$.

Let us denote by ${\rm ch}(V^{\tilde\la}_{m|n})$ the formal
character of the irreducible $gl_{m|n}$-module. Theorem
\ref{hw-glmn} leads to
\begin{equation*}
{\rm ch}(\C[\x,\etabf]) = \sum_{\la} {\rm ch}(V^{\la}_d) {\rm
ch}(V^{\tilde\la}_{m|n}),
\end{equation*}
where, we recall that,  ${\rm ch}(V^{\la}_d)=s_{\la}(\bar{\bf
z})$. By using \lemref{independence} we obtain the following
well-known result \cite{BR}.
\begin{thm}\label{characterthm} For each partition $\la$ of length
$d$ subject to the condition $\la_{m+1}\le n$,
\begin{equation*}
{\rm ch}V_{m|n}^{\widetilde\la}=HS_{\la}(\bar \x;\bar \etabf),
\end{equation*}
where $\widetilde\la$ is defined by \eqref{wtglmn}.
\end{thm}

Keep the notations of this subsection but replace $m$ by $p$ and
$n$ by $q$. Let $\bar\x^{-1}$ $=$ $\{\bar x_1^{-1},\bar
x_2^{-1},\cdots ,\bar x_p^{-1}\}$, $\bar \etabf^{-1}=\{\bar
\eta_1^{-1},\bar \eta_2^{-1},\cdots ,\bar \eta_q^{-1}\}$ and $\bar
{\bf z}^{-1}=\{\bar z_1^{-1},\bar z_2^{-1},\cdots,\bar
z_d^{-1}\}$. Using \eqref{barphid} and \eqref{barphi}, we can
easily compute the formal character of the $gl_d\times
gl_{p|q}$-module $\C[\y,\zetabf]$ with respect to the Cartan
subalgebra $\sum_{i=1}^{p+q} \C E_{ii} \oplus\sum_{k=1}^d\C
e_{kk}$. We have
\begin{equation}
{\rm ch}(\C[\y,\zetabf]) =\left(\bar x_1\cdots\bar x_p\right)^{-d}
\left(\bar\eta_1\cdots\bar\eta_q\right)^d \prod_{i,j,k}(1-\bar
x_i^{-1}\bar z_k^{-1})^{-1}(1+\bar \eta_j^{-1}\bar z_k^{-1}),
\end{equation}
where $1\le i\le p$, \  $1\le j\le q$, and $1\le k\le d$. By
\propref{combid*},
\begin{equation}\label{character2}
{\rm ch}(\C[\y,\zetabf])=\left(\bar x_1\cdots\bar x_p\right)^{-d}
\left(\bar\eta_1\cdots\bar\eta_q\right)^d \sum_{\la}HS_{\la}(\bar
\x^{-1};\bar \etabf^{-1})s_{\la}(\bar{\bf z}^{-1}),
\end{equation}
where $\la$ is summed over all partitions with length $d$ subject
to $\la_{p+1}\le q$.

Note that ${\rm ch}(V_d^{\la^*})=s_{\la}(\bar{\bf z}^{-1})$. Thus
the following theorem is a consequence of \thmref{hw*} by using
\lemref{independence} and equation \eqref{character2}.
\begin{thm}\label{characterthm*}
For each partition $\la$ of length $d$ subject to the condition
$\la_{p+1}\le q$,
\begin{equation*}
{\rm ch}V_{p|q}^{-d {\bf 1}+\widehat{\la^*}}= \left(\bar
x_1\cdots\bar x_p\right)^{-d}
\left(\bar\eta_1\cdots\bar\eta_q\right)^d HS_{\la}(\bar
\x^{-1};\bar \etabf^{-1}),
\end{equation*}
where $\widehat{\la^*}$ is as given in \eqref{glpqwts}.
\end{thm}

\subsection{Character formulas for unitarizable
$gl_{m+p|n+q}$-modules}

We keep the notations $z_i$, $1\le i\le d$, and ${\bf z}$, \ ${\bf
z}^{-1}$ from the last subsection. Let $e$ be the formal
indeterminate as before. For $1\le r\le p,$ \ $1\le s\le q,$ \
$1\le i\le m,$ and $1\le j\le n$, we define
\begin{eqnarray*}
\bar y_r=e^{\widehat\epsilon_r},&\quad&
\bar\zeta_s=e^{\widehat\epsilon_{s+p}},   \\
\bar x_i=e^{\widehat\epsilon_{i+p+q}}, &\quad&
\bar\eta_j=e^{\widehat\epsilon_{j+m+p+q}},
\end{eqnarray*}
where, we recall that, $\widehat\epsilon_A\in\fh_{m+p|n+q}^*$,
$A\in{\bf I}$, are defined by $\widehat\epsilon_A(E^B_B)=\delta_{A
B}$, $A, B\in{\bf I}$. Set $\bar\x=\{\bar x_1,\bar x_2,\cdots
,\bar x_m\}$, $\bar \etabf=\{\bar \eta_1,\bar \eta_2,\cdots ,\bar
\eta_n\}$, $\bar\y^{-1}=\{\bar y_1^{-1},\bar y_2^{-1},\cdots ,\bar
y_p^{-1}\}$ and $\bar \zetabf^{-1}=\{\bar \zeta_1^{-1},\bar
\zeta_2^{-1},\cdots ,\bar \zeta_q^{-1}\}$.

We wish to compute the formal characters ${\rm
ch}(W_{m+p|n+q}^{\La(\la)})$  with respect to the Cartan
subalgebra $\fh_{m+p|n+q}=\sum_{A=1}^{m+n+p+q} \C E^A_A$ for the
unitarizable irreducible $gl_{m+p|n+q}$-modules
$W_{m+p|n+q}^{\La}$  appearing in the decomposition of $\C[\x, \y,
\etabf, \zetabf]$.

\begin{thm}\label{chWLa}
For each generalized partition $\la$ of length $d$
subject to the conditions $\la_{m+1}\le n$ and $\la_{d-p}\ge -q$,
\begin{equation*}
{\rm ch}(W_{m+p|n+q}^{\La(\la)})=(\bar y_1\bar y_2\cdots\bar
y_p)^{-d}(\bar \zeta_1\bar \zeta_2\cdots\bar \zeta_q)^{d}
\sum_{\mu,\nu} C^\la_{\mu \nu^*} HS_{\mu}(\bar \x;\bar
\etabf)HS_{\nu}(\bar \y^{-1};\bar \zetabf^{-1}),
\end{equation*}
where $\mu$ and $\nu$ are summed over all partitions of length $d$
subject to the conditions $\mu_{m+1}\le n$ and $\nu_{p+1}\le q$
respectively. The $C^\la_{\mu \nu^*}$ are the
Littlewood-Richardson coefficients.
\end{thm}
\begin{proof}
Consider the restriction of $W_{m+p|n+q}^{\La(\la)}$ to a module of the
subalgebra $gl_d\times gl_{m|n}\times gl_{p|q}$. Its formal
character with respect to the Cartan subalgebra
$\fh_d\times\fh_{p|q}\times\fh_{m|n}$ coincides with ${\rm
ch}(W_{m+p|n+q}^{\La(\la)})$. Therefore by \thmref{branchrule} we have
\begin{equation}
{\rm ch}(W_{m+p|n+q}^{\La(\la)})= \sum_{\mu,\nu}C^\la_{\mu \nu^*}
{\rm ch}(V_{m|n}^{\tilde\mu}){\rm ch}(V_{p|q}^{-d {\bf 1}+
\widehat{\nu^*}})
\end{equation}
where $\mu$ and $\nu$ are summed over all partitions of length $d$
subject to the conditions $\mu_{m+1}\le n$ and $\nu_{p+1}\le q$
respectively.  Using \thmref{characterthm} and
\thmref{characterthm*} in this equation we immediately arrive at
the claimed result.
\end{proof}

\section{Tensor Product Decomposition of unitarizable irreducible
$gl_{m+p|n+q}$-modules}\label{tensorproduct}

As another application of \thmref{hw+*}, we shall compute the
tensor product decomposition
\begin{equation}\label{tensor}
W_{m+p|n+q}^{\La(\mu)}\otimes
W_{m+p|n+q}^{\La(\nu)}\cong\sum_{\la}a_{\mu\nu}^\la
W_{m+p|n+q}^{\La(\la)},
\end{equation}
where $\mu$ and $\nu$ are generalized partitions of length $l$ and
$r$, respectively, satisfying in addition the conditions
$\mu_{m+1}\le n$, $\nu_{m+1}\le n$, $\mu_{l-p}\ge -q$ and
$\nu_{r-p}\ge -q$. It will follow easily from our discussion that
the summation $\la$ in \eqnref{tensor} is over all generalized
partitions of length $l+r$ and satisfying $\la_{m+1}\le n$ and
$\la_{l+r-p}\ge -q$. We will compute the coefficients
$a_{\mu\nu}^\la$ in terms of the usual Littlewood-Richardson
coefficients (see e.g.~\cite{M}).

We have by \thmref{hw+*} for $d=l$, $r$ and $l+r$ respectively
\begin{equation*}
S(\C^l\otimes\C^{m|n}\oplus\C^{l*}\otimes\C^{p|q*})\cong\sum_{\mu}V^\mu_l\otimes
W^{\La(\mu)}_{m+p|n+q},
\end{equation*}
\begin{equation*}
S(\C^r\otimes\C^{m|n}\oplus\C^{r*}\otimes\C^{p|q*})\cong\sum_{\nu}V^\nu_r\otimes
W^{\La(\nu)}_{m+p|n+q},
\end{equation*}
and
\begin{equation*}
S(\C^{l+r}\otimes\C^{m|n}\oplus\C^{l+r*}\otimes\C^{p|q*})\cong\sum_{\la}V^\la_{l+r}\otimes
W^{\La(\la)}_{m+p|n+q},
\end{equation*}
where $\mu$, $\nu$ and $\la$ are generalized partitions satisfying
the corresponding conditions described above.   The isomorphism
$S(\C^l\otimes\C^{m|n}\oplus\C^{l*}\otimes\C^{p|q*})\otimes
S(\C^r\otimes\C^{m|n}\oplus\C^{r*}\otimes\C^{p|q*})\cong
S(\C^{l+r}\otimes\C^{m|n}\oplus\C^{l+r*}\otimes\C^{p|q*})$ gives
rise to
\begin{equation}\label{aux61}
\sum_{\mu,\nu}V^\mu_l\otimes V^\nu_r\otimes W^{\La(\mu)}_{m+p|n+q}
\otimes W^{\La(\nu)}_{m+p|n+q}\cong \sum_{\la}V^\la_{l+r}\otimes
W^{\La(\la)}_{m+p|n+q}.
\end{equation}

Now suppose that $V_{l+r}^\la$, when regarded as a $gl_l\times
gl_r$-module via the obvious embedding of $gl_l\times gl_r$ into
$gl_{l+r}$, decomposes as
\begin{equation*}
V_{l+r}^\la\cong \sum_{\mu,\nu}b_\la^{\mu\nu}V_l^\mu\otimes
V_r^\nu.
\end{equation*}
This together with \eqnref{tensor} and \eqnref{aux61} give
\begin{equation}\label{aux62}
a_{\mu\nu}^\la=b_\la^{\mu\nu}.
\end{equation}

The duality between the branching coefficients and tensor products
of a general dual pair is well-known \cite{H2}.  We recall that in
\eqnref{aux62} $\mu$, $\nu$ and $\la$ are generalized partitions
subject to the appropriate constraints.

Now \thmref{hw-glmn} with $n=0$ combined with an analogous
argument as the one given above imply that
\begin{equation}\label{aux63}
C_{\mu\nu}^\la=b_\la^{\mu\nu},
\end{equation}
where here $\mu$, $\nu$, $\la$ are partitions of appropriate
lengths and the $C_{\mu\nu}^\la$'s are the usual
Littlewood-Richardson coefficients.

Now for generalized partitions $\mu$, $\nu$ and $\la$ subject to
appropriate constraints the decomposition $V_{l+r}^\la\cong
\sum_{\mu,\nu}b_\la^{\mu\nu}V_l^\nu\otimes V_r^\mu$ implies that
$V_{l+r}^{\la+d{\bf 1}_{l+r}}\cong
\sum_{\mu,\nu}b_\la^{\mu\nu}V_l^{\mu+d{\bf 1}_l}\otimes
V_r^{\nu+d{\bf 1}_r}$, where ${\bf 1}_k$ denotes the $k$-tuple
$(1,1,\cdots,1)$ regarded as a partition.  Hence
$b_\la^{\mu\nu}=b_{\la+d{\bf 1}_{l+r}}^{\mu+d{\bf 1}_l,\nu+d{\bf
1}_r}$. Now if we choose a non-negative integer $d$ so that
${\la+d{\bf 1}_{l+r}}$ is a partition, then $b_{\la+d{\bf
1}_{l+r}}^{\mu+d{\bf 1}_l,\nu+d{\bf 1}_r}=C^{\la+d{\bf
1}_{l+r}}_{\mu+d{\bf 1}_l,\nu+d{\bf 1}_r}$ and hence by
\eqnref{aux62} and \eqnref{aux63}
\begin{equation*}
a^\la_{\mu\nu}=C^{\la+d{\bf 1}_{l+r}}_{\mu+d{\bf 1}_l,\nu+d{\bf
1}_r}.
\end{equation*}

From our discussion above we arrive at the following theorem.

\begin{thm}\label{tensor1} Let $\mu$ and $\nu$ be generalized
partitions of length $l$
and $r$, respectively, satisfying in addition the conditions
$\mu_{m+1}\le n$, $\nu_{m+1}\le n$, $\mu_{l-p}\ge -q$ and
$\nu_{r-p}\ge -q$.  Let $W^{\La(\mu)}_{m+p|n+q}$ and
$W^{\La(\nu)}_{m+p|n+q}$ be the corresponding unitarizable
$gl_{m+p|n+q}$-modules. We have the following decomposition of
$W_{m+p|n+q}^{\La(\mu)}\otimes W_{m+p|n+q}^{\La(\nu)}$ into
irreducible $gl_{m+p|n+q}$-modules:
\begin{equation*}
W_{m+p|n+q}^{\La(\mu)}\otimes
W_{m+p|n+q}^{\La(\nu)}\cong\sum_{(\la,d)}C_{\mu+d{\bf
1}_l,\nu+d{\bf 1}_r}^\la W_{m+p|n+q}^{\La(\la-d{\bf 1}_{l+r})},
\end{equation*}
where the summation above is over all pairs $(\la,d)$ subject to
the following four conditions:
\begin{itemize}
\item[(i)] $\la$ is a partition of length $l+r$ and $d$ a
non-negative integer. \item[(ii)] $(\la-d{\bf 1}_{l+r})_{m+1}\le
n$ and $(\la-d{\bf 1}_{l+r})_{l+r-p}\ge -q$. \item[(iii)]
$\mu+d{\bf 1}_l$ and $\nu+d{\bf 1}_r$ are partitions. \item[(iv)]
If $d>0$, then $\la$ is a partition with $\la_{l+r}=0$.
\end{itemize}
Here the coefficients $C_{\mu+d{\bf 1}_l,\nu+d{\bf 1}_r}^\la$ are
determined by the tensor product decomposition of $gl_{k}$-modules
$V^{\mu+d{\bf 1}_l}_{k}\otimes V^{{\nu+d{\bf
1}_r}}_{k}\cong\sum_{\la}C_{\mu+d{\bf 1}_l,\nu+d{\bf 1}_r}^\la
V^\la_{k}$, where $k= l+r$.
\end{thm}

\begin{rem}
In the above theorem the coefficients $C_{\mu+d{\bf 1}_l,\nu+d{\bf
1}_r}^\la$ are the usual Littlewood-Richardson coefficients
associated to partitions, and hence can be computed via the
Littlewood-Richardson rule.
\end{rem}

\begin{rem}
We finally note the remarkable similarity between the irreducible
representations of the so-called super $W_{1+\infty}$, that is the
Lie superalgebra of differential operators on the circle with
$N=1$ extended symmetry, that appear in the decomposition of
tensor powers of its natural representation on the
infinite-dimensional Fock space \cite{CL}, and the irreducible
unitarizable representations of $gl_{m+p|n+q}$ of this paper.
Indeed the characters and the tensor product decomposition are
virtually identical modulo some modification necessitated by the
infinite-dimensional nature of the super $W_{1+\infty}$. This
similarity can be explained by the existence a Howe duality
between the super $W_{1+\infty}$ and $gl_d$ on the $d$-th tensor
power of the Fock space generated by infinitely many fermionic and
bosonic quantum oscillators \cite{CW3}.
\end{rem}


\end{document}